\documentclass[final,3p,times]{elsarticle}

 \usepackage{graphics}
 \usepackage{graphicx}
 \usepackage{epsfig}
\usepackage{caption}
\usepackage{subcaption}
\usepackage{varwidth}
\usepackage{multirow}
\usepackage{hhline}
\usepackage{wrapfig}
\usepackage[boxruled]{algorithm2e}
\usepackage{float}
\usepackage{diagbox}

\usepackage{geometry}
\usepackage[table]{xcolor}
\usepackage{booktabs,tabularx}
\renewcommand\toprule{\specialrule{1pt}{1pt}{0pt}\rowcolor{gray!10}}
\renewcommand\midrule{\specialrule{0.4pt}{0pt}{0pt}}

\DeclareCaptionFormat{myformat}{\fontsize{10}{15}\selectfont#1#2#3}
\usepackage{color}
\usepackage{amsmath}
\usepackage{amssymb}
\usepackage{amsthm}
\newtheorem{thm}{Theorem}

\newtheorem*{remark}{Remark}

\newtheorem{example}{\textbf{Example}}

\newcommand{\s}{{\mathbf s}}

\journal{NA}

\begin{document}

\begin{frontmatter}



\title{Hierarchical randomized neural networks for efficiently and accurately finding multiple solutions of semilinear elliptic equations}
\author{Jingying Hu\fnref{label1,label4}}
\author{Jiao Li\fnref{label3}}
\author{Yufeng Xu\fnref{label2}}
\author{Jinyong Ying\fnref{label2}\corref{cor1}}

\fntext[label1]{Shenzhen Institutes of Advanced Technology, Chinese Academy of Sciences, Shenzhen, Guangdong 518055, China}
\fntext[label4]{University of Chinese Academy of Sciences, Beijing, 101408, China}
\fntext[label3]{School of Mathematics and Statistics, Hunan Provincial Key Laboratory of Mathematical Modeling and Analysis in Engineering, Changsha University of Science and Technology, Changsha, Hunan 410114, China}
\fntext[label2]{School of Mathematics and Statistics, HNP-LAMA, Central South University, Changsha, Hunan 410083, China}
 \cortext[cor1]{Corresponding author, email address: ying48882345@gmail.com}
\begin{abstract}
Being able to find multiple solutions of a semilinear elliptic equation is a fundamental problem in scientific computing. In this paper, fully taking the advantages of randomized neural networks, we propose a new efficient and accurate method for identifying potentially as many solutions as possible. The new method is constructed from three aspects--initialization, optimization, a technique for acceleration and improving accuracy, two of which are innovative developments. Firstly, based on the basis functions of a randomized  network and the sampling points, a characteristic matrix is introduced so that its singular value decomposition can be used to determine weakly variant search directions, followed by a newly developed algorithm to generate  effective initial candidates. Then the good deflated Gauss-Newton method [see {\em A. Riley, et al., SIAM Journal on Scientific Computing 47~(5) (2025) A2963-A2987}] using an exponential deflation operator is adapted as the optimization to search for all multiple solutions. Local error estimate of a particular solution under some proper condition is given in terms of the neurons in the hidden layer and the sampling points for this newly designed method. Furthermore, to accelerate the computation and improve the solution accuracy, a new multilevel method based on hierarchical networks is designed.  Without using any {\em a priori} information, several numerical tests show together with the multilevel method, the new method can not only accurately find multiple solutions with relative errors reaching $10^{-16}$ only in less than one-tenth computational cost for the best case, but also preserve the inherent symmetries of multiple solutions. 
\end{abstract}


\end{frontmatter}


\section{Introduction}
Nonlinear elliptic equations arise in a wide range of scientific and engineering disciplines,
including quantum physics, solid mechanics, and fluid dynamics, where they model steady-state or
equilibrium phenomena. 
The properties of their solutions 
are closely tied to the physical behavior of the underlying systems. Although one-dimensional
problems have been thoroughly studied using analytical tools such as Galerkin projections and
energy estimates, the analysis of multidimensional or highly nonlinear cases remains challenging.
Consequently, developing efficient and reliable numerical methods for approximating solutions to
nonlinear elliptic problems is not only of theoretical interest but also a critical demand in
scientific computing and practical engineering applications.

This is the first attempt of our serial work on finding multiple solutions of nonlinear elliptic
equation. Let us begin with the nonlinear elliptic equation with polynomial nonlinearity.
In this paper, we are concerned with the following nonlinear elliptic equation
\begin{equation}\label{Eq0101}
\left\{\begin{array}{cl}
  -\Delta{u}(x) = \lambda f(x,u),~~&~~x\in\Omega\subset\mathbb{R}^n,\\
  u(x) = 0,~~&~~x\in\partial\Omega,
\end{array}\right.
\end{equation}
where $\Omega \in \mathbb{R}^n$ is some bounded domain, $f(x,u)$ is
a polynomial nonlinearity with constant or variable coefficients, $\lambda>0$ is a real
positive constant. This is a typical mathematical model contains the Lane-Emden equation, Allen-Cahn
equation, nonlinear eigenvalue problem and elliptic Lotka-Volterra system as special cases.
 The study of existence
and multiplicity of solutions of model \eqref{Eq0101} by using variational method and other
nonlinear analysis techniques has received considerable attention, given that $f$ satisfies
some necessary increasement restrictions. For a comprehensive review, we refer
the readers to an interesting small monograph by David Costa \cite{DGC2007}.  
It illustrates many examples with multiple solutions, which also inspire us to further ask:
\emph{How many classical solutions model \eqref{Eq0101} will have exactly? If the number is finite but not unique,
or infinite, can we find all of these classical solution numerically?}
While theoretical results are usually discussed in a functional space with general dimensionality $n$, for the convenience
of simulations, we shall focus on $1 \leq n\leq 3$.

For the model problem, comparing to the associate linear problem $-\Delta u(x) = \lambda u(x)$, $x\in\Omega\subset
\mathbb{R}^n$ with $u(x)=0$, $x\in\partial\Omega$, whose distribution of eigenvalues and
behavior of eigenfunction for regular domain are clear, the semilinear elliptic equation with
general polynomial nonlinearity is not well-studied yet as it involves more complicated dynamics, especially when
$f(x,u)$ is some particular cubic polynomial of $u$. For example, 
in \cite{AC2000}, the bounded solutions of semilinear elliptic equations $\Delta u - F'(u)=0$ in
the whole space $\mathbb{R}^n$ is studied. The investigation of existence of solutions to this problem
with $n\leq 8$ drives a famous conjecture by De Giorgi.
About ten year later, Pino, Kowalczyk and Wei find a counterexample to De Giorgi's conjecture
for dimension $N\geq 9$ \cite{PKW2011}.
In \cite{Tang2025}, the author proves that for nonlinearity $f(u)=|u|^{p-1}u-u$,
$1<p<\frac{n+2}{n-2}$, the elliptic equation $-\Delta u = f(u)$ in $\mathbb{R}^n$ $(n\geq 3)$
with $u\in H^1(\mathbb{R}^n)$ and $u\nequiv 0$ has a unique solution $u=u(|x|)$. 
These celebrate progresses motivate us to consider what property the model \eqref{Eq0101} will have theoretically if it defines on a smooth bounded domain $\Omega$ with
general polynomial nonlinearity.

From a numerical viewpoint, solving nonlinear elliptic equations still runs into several challenges.
First of all, the inherent nonlinearity of these problems often leads to non-convex energy functionals,
multiple solutions, and bifurcation structures, rendering classical linear numerical methods
ineffective.
Secondly, two-dimensional and even high-dimensional versions of these equations yield
large-scale discrete systems that require advanced computational strategies to maintain stability
and efficiency.
Thirdly, an additional difficulty arises from multiscale behavior, boundary layers, and
singularities, necessitating numerical schemes with strong adaptivity and robustness.
Classic methods such as finite element method, Galerkin-based methods together with Newton-type techniques, homotopy continuation, and
multilevel technique \cite{AB77,GT00} have shown some promising abilities to resolve some of the mentioned issues.  Recently, due to the success of neural networks and the improvements for solving various problems \cite{raissi2019physics,ying2024improved}, physics-informed neural networks have also been used to discover multiple solutions\cite{zou2025learning,huang2022hompinns,disaro2027deflation}, where neural networks are mainly used to find good initializations and then classic methods are again used to maintain high accuracy. Overall, for arbitrary semilinear problems, efficiently and accurately finding all solutions of a general nonlinear elliptic equation is still far to be done.



Therefore, in this paper, we mainly focus on the investigation of the numerical method for the model problem with multiple solutions. Instead of using classic methods, we fully take the advantages of randomized neural networks and propose a new efficient and accurate method for identifying potentially as many solutions as possible without {\em a priori} information, although some analysis results can speed up finding these solutions. The new method is constructed from three different aspects--initializations, optimizations, and improvements in terms of efficiency and accuracy. 
And main contributions of this paper are as follows:
\begin{enumerate}
\item Based on randomized neural networks and sampling points, the characteristic matrix and the weakly variant search directions are introduced. Meanwhile, the singular value decomposition of the characteristic equation is used to generate effective candidates at initialization step, instead of generating initialization randomly, which is very low efficiency;
\item New initialization algorithm together with the recently developed good deflated Gauss-Newton method is combined to form a new algorithm. Under some proper conditions for some solutions, the local error estimate of the randomized neural network is presented in terms of numbers of both neurons in the hidden layer and the sampling points;
\item Multilevel method, which is like the full multigrid method in some sense,  is further proposed to accelerate the computation and improve the solution accuracy;
\item Various examples with or without exact solutions are used to demonstrate the effectiveness of the proposed method, reaching the relative error $10^{-16}$ for the best case while reducing the computational cost to only less than $8\%$. Meanwhile, numerical tests show that the proposed method is able to preserve the inherent symmetries of multiple solutions, which is a very important property for many applications.
\end{enumerate}

The rest paper is organized as follows. The second section is the introduction of the considered model, some analysis could be used to prove existence of different solutions, and the difficulties of identifying number of solutions by pure analysis tools. The third section is devoted to the numerical method to identify multiple solutions, including randomized neural networks, new algorithm for initialization, the optimization algorithm, and the local error estimate for a solution satisfying some conditions. Then a new improved method based  on hierarchical networks, the multilevel method, is proposed in Section 4, followed by several numerical tests from one-dimensional cases to three dimensional ones in Section 5, showing the effective and improvements of the proposed method. Lastly a conclusion together with some future work is given.

\section{Model description}
In this paper, we shall recall the elementary analysis of nonlinear elliptic equations.
Firstly, we will show the variational property of model problem and introduce some existence
results of weak solutions. Then we present the mathematical analysis of multiplicity
for one-dimensional nonlinear elliptic equation.

\subsection{Variational formulation}
Let $\Omega\subset\mathbb{R}^n$ be a bounded open set with $C^1$-boundary, and we 
consider the Dirichlet problem \eqref{Eq0101}
with $f:\Omega\times\mathbb{R}\to\mathbb{R}$ being a Carath\'eodory function. For $n>2$,
one usually writes $2^*=\frac{2n}{n-2}$, while for $n\leq 2$, it becomes $2^*=+\infty$.
We denote $H=H_0^1(\Omega)$ equipped with the Hilbert norm
$\|u\|=\left(\int_\Omega |\nabla u|^2 dx\right)^{1/2}$
and the corresponding inner product. By the Poincar\'e inequality,
this norm is equivalent to the standard $H^1(\Omega)$-norm on $H_0^1(\Omega)$.
The continuous embedding $H\hookrightarrow L^q(\Omega)$, $1\leq q\leq 2^*$ and
compact embedding $H\hookrightarrow L^q(\Omega)$, $1\leq q<2^*$ are useful in the discussion later.
The dual space of \(H\) is denoted by \(H^{-1}(\Omega)\).

A function \(u\in H\) is called a weak solution of the problem \eqref{Eq0101} if
\[
\int_\Omega \nabla u\cdot\nabla v dx
=
\lambda\int_\Omega f(x,u)v\,dx
\qquad\text{for every }v\in H.
\]
Define $F(x,s)=\int_0^s f(x,t) dt$ and the energy functional
\[
I_\lambda(u)=\frac12\|u\|^2-\lambda\int_\Omega F(x,u) dx,
\qquad u\in H.
\]
Under the growth assumptions below, \(I_\lambda\in C^1(H,\mathbb{R})\) and
\[
I_\lambda'(u)v
=
\int_\Omega \nabla u\cdot\nabla v dx
-
\lambda\int_\Omega f(x,u)v dx,
\qquad v\in H.
\]
Thus weak solutions are precisely the critical points of \(I_\lambda\). Firstly, we have the following result, proof of which is omitted as it is standard. 
\begin{thm}
Assume that $f:\Omega\times\mathbb{R}\to\mathbb{R}$ is Carath\'eodory, and there exists a constant $M>0$, such that
\[
|f(x,s)|\le M \qquad\text{for a.e. }x\in\Omega\text{ and all }s\in\mathbb{R}.
\]
Then for every $\lambda>0$, the model problem \eqref{Eq0101}
has at least one weak solution $u\in H_0^1(\Omega)$. 
Moreover, such a solution is the global minimizer of \(I_\lambda\) on \(H_0^1(\Omega)\).
\end{thm}

It is easy to see that $|f(x,s)|\leq M$ is too restrictive for polynomial nonlinearity,
although it holds whenever $s$ is finite and pre-given. We shall consider the existence of solution
to the model problem in weaker assumptions.

\begin{thm}
Assume that $f:\Omega\times\mathbb{R}\to\mathbb{R}$ is Carath\'eodory and 
satisfies the following conditions:
\begin{itemize}
\item[(G1)] There exist constants $C>0$ and $1<p<2^*-1$, such that
\[
|f(x,s)|\leq C(1+|s|^p)
\qquad\text{for a.e. }x\in\Omega\text{ and all }s\in\mathbb{R}.
\]
If $n\leq 2$, any finite $p\geq 1$ is allowed.

\item[(G2)] Uniformly in $x\in\Omega$,
\[
f(x,s)=o(|s|)
\qquad\text{as }s\to0.
\]

\item[(G3)] There exist constants $\theta>2$ and $R>0$, such that
\[
0<\theta F(x,s)\leq s f(x,s)
\qquad\text{for a.e. }x\in\Omega\text{ and all }|s|\geq R.
\]
\end{itemize}
Then for every $\lambda>0$, the model problem \eqref{Eq0101} 
has at least one nontrivial weak solution $u\in H_0^1(\Omega)$.
\end{thm}

\begin{proof}
By assumption (G1), the functional $I_\lambda(u)$
is well defined and of class $C^1$ in $H$, and
\[
I_\lambda'(u)v
=
\int_\Omega \nabla u\cdot\nabla v dx
-
\lambda\int_\Omega f(x,u)v dx,
\qquad v\in H.
\]
We first verify the famous Palais-Smale condition. Let $\{u_n\}\subset H$ be a sequence, 
such that
\[
I_\lambda(u_n)\to c
\qquad\text{and}\qquad
I_\lambda'(u_n)\to0\quad\text{in }H^{-1}(\Omega).
\]
Using assumption (G3), we compute
\[
\theta I_\lambda(u_n)-I_\lambda'(u_n)u_n
=
\left(\frac{\theta}{2}-1\right)\|u_n\|^2
-
\lambda\int_\Omega
\left(\theta F(x,u_n)-f(x,u_n)u_n\right) dx.
\]
It is easy to have that for $|s|\geq R$, condition (G3) implies $
\theta F(x,s) - s f(x,s) \le 0$. 
Similarly, for $|s|\leq R$, the growth condition (G1) gives that $
|\theta F(x,s)-s f(x,s)| \leq C_R$ 
for some constant $C_R>0$ depending only on $R$. 
Hence
\[
\int_\Omega
\left(\theta F(x,u_n)-f(x,u_n)u_n\right) dx
\leq C_R|\Omega|, 
\]
from which it implies
\[
\left(\frac{\theta}{2}-1\right)\|u_n\|^2
\leq
\theta I_\lambda(u_n)-I_\lambda'(u_n)u_n
+
\lambda C_R|\Omega|.
\]
Since $I_\lambda(u_n)$ is bounded and $I_\lambda'(u_n) \to 0$ in $H^{-1}(\Omega)$,
$
\left(\frac{\theta}{2}-1\right)\|u_n\|^2
\leq
C(1+\|u_n\|).
$
Because $\theta > 2$, it follows that $\{u_n\}$ is bounded in $H$. Thus, after passing to a subsequence, $
u_n\rightharpoonup u$ weakly in $H$.
By the compact embedding $H\hookrightarrow L^{p+1}(\Omega)$, we have
\[
u_n\to u\quad\text{in }L^{p+1}(\Omega).
\]
From the growth condition (G1), the sequence $\{f(x,u_n)\}$ 
is bounded in $L^{(p+1)/p}(\Omega)$. 
Hence, by H\"older's inequality,
\[
\left|\int_\Omega f(x,u_n)(u_n-u) dx\right|
\leq
\|f(x,u_n)\|_{L^{(p+1)/p}(\Omega)}
\|u_n-u\|_{L^{p+1}(\Omega)}
\to0.
\]
Since $I_\lambda'(u_n) \to 0$ in $H^{-1}(\Omega)$, we have $I_\lambda'(u_n)(u_n-u)  \to  0$.
Therefore, $\int_\Omega \nabla u_n\cdot\nabla(u_n-u) dx  \to 0$.
That is,
$
\|u_n\|^2-\int_\Omega \nabla u_n\cdot\nabla u dx \to 0.
$
But weak convergence gives
$
\int_\Omega \nabla u_n\cdot\nabla u dx \to \|u\|^2.
$
Hence we know
$
\|u_n\|^2\to\|u\|^2.
$
Together with weak convergence, this implies
\[
u_n\to u  \quad  \text{strongly in}~~H.
\]
Consequently $I_\lambda'(u)=0$. 
Thus $I_\lambda$ satisfies the Palais-Smale condition.

Next we prove the {\it Mountain Pass Geometry}. 
From conditions (G2) and (G1), for every $\varepsilon>0$, there exists a constant 
$C_\varepsilon>0$, such that
\[
|F(x,s)|\leq \frac{\varepsilon}{2}|s|^2 + C_\varepsilon|s|^{p+1}.
\]
Therefore
\[
I_\lambda(u) \geq  \frac12\|u\|^2
-
\lambda\frac{\varepsilon}{2}\|u\|_{L^2(\Omega)}^2
-
\lambda C_\varepsilon\|u\|_{L^{p+1}(\Omega)}^{p+1}.
\]
Using the continuous embedding results 
$H\hookrightarrow L^2(\Omega)$ and $H\hookrightarrow L^{p+1}(\Omega)$, 
we obtain
\[
I_\lambda(u) \geq
\frac12\left(1-\lambda\varepsilon C_P^2\right)\|u\|^2
-
C_\varepsilon'\|u\|^{p+1}.
\]
We may choose $\varepsilon>0$ so small that
$
1-\lambda\varepsilon C_P^2 > 0.
$
Since $p+1>2$, it follows that for sufficiently small $\rho>0$,
\[
I_\lambda(u)  \geq \alpha>0
\qquad\text{whenever}~~\|u\| = \rho.
\]
Also $I_\lambda(0)=0$. Thus $0$ is a local minimum of $I_\lambda$.

On the other hand, condition (G3) implies the existence of constants $a>0$ and $b\geq 0$, 
such that
\[
F(x,s)\geq a|s|^\theta - b
\qquad\text{for a.e. }x \in \Omega\text{ and all }  s\in\mathbb{R}.
\]
We choose $\varphi\in C_c^\infty(\Omega)$ with $\varphi \geq 0$ and 
$\varphi\not\equiv 0$. Then for $t>0$,
\[
I_\lambda(t\varphi) \leq \frac{t^2}{2}\|\varphi\|^2
-
\lambda a t^\theta\int_\Omega |\varphi|^\theta dx
+
\lambda b|\Omega|.
\]
Since $\theta>2$, the right-hand side tends to $-\infty$ as $t \to +\infty$. 
Hence for sufficiently large $t>0$,
$
I_\lambda(t\varphi) < 0.
$
Set $e=t\varphi$ for such a given $t$. Then $I_\lambda(e)<0$. Therefore $I_\lambda$ satisfies the hypotheses of the mountain pass theorem (see \cite{DGC2007}). 
Hence $I_\lambda$ has a critical point $u\in H$ with
$
I_\lambda(u) \geq \alpha >0.
$
In particular, $u \neq 0$. 
Since every critical point of $I_\lambda$ is a weak solution of the model problem, 
$u$ is a nontrivial weak solution.
\end{proof}

Based on Theorem 1 and Theorem 2, it is still difficult to know the number of solutions 
to the model problem in general. Next, we consider a one-dimensional model problem with
polynomial term $f(u)$ at most quintic, to show how complex this task will be.

\subsection{Mathematical analysis for one-dimensional model problem}
In this section, we discuss some theoretical results of nonlinear elliptic equation with
polynomial nonlinearity for one-dimensional case. Let us begin with
\begin{equation}\label{Eq0203}
\left\{\begin{array}{cl}
  -{u}''(x) =  \alpha_1 u + \alpha_2 u^2 + \alpha_3 u^3 + \alpha_4 u^4 + \alpha_5 u^5,
  ~~&~~x\in(0,1),\\
  u(0) = u(1) = 0,
\end{array}\right.
\end{equation}
where $\Omega = (0,1)$, and $\alpha_k$ $(k=1,2,3,4,5)$ are real constants.  Let $v(x)=u'(x)$, and we have
\begin{equation}\label{Eq0204}
\left\{\begin{array}{ll}
  u'(x) = v(x), \\
  v'(x) =  - \alpha_1 u - \alpha_2 u^2 - \alpha_3 u^3 - \alpha_4 u^4 - \alpha_5 u^5
\end{array}\right.
\end{equation}
with $u(0)=0$ and $v(0)=\eta$, where $\eta$ is a constant determined later. Direct computation shows
that
\[
\frac{du}{dv} = -\frac{v(x)}{\alpha_1 u + \alpha_2 u^2 + \alpha_3 u^3 + \alpha_4 u^4 + \alpha_5 u^5},
\]
which implies
\[
\frac{1}{2}[v(x)]^2 = \frac{1}{2}\eta^2 - \int^{u(x)}_{u(0)}
\left(\alpha_1 u + \alpha_2 u^2 + \alpha_3 u^3 + \alpha_4 u^4 + \alpha_5 u^5\right)du.
\]
Hence, we arrive at
\begin{equation}\label{Eq0205}
\frac{du}{dx} = \pm\sqrt{\eta^2 - \left(\alpha_1 u^2 + \frac{2\alpha_2}{3} u^3
+ \frac{\alpha_3}{2} u^4 + \frac{2\alpha_4}{5} u^5 + \frac{\alpha_5}{3} u^6\right)}.
\end{equation}
Integrating from $[0,x)$ with respect to variable $x$, we get
\begin{equation}\label{Eq0206}
  \int^{u(x)}_{u(0)} \frac{du}{\pm\sqrt{\eta^2 - \left(\alpha_1 u^2 + \frac{2\alpha_2}{3} u^3
+ \frac{\alpha_3}{2} u^4 + \frac{2\alpha_4}{5} u^5 + \frac{\alpha_5}{3} u^6\right)}} = x.
\end{equation}
Therefore, we may see that the key step to search an elementary solution to model \eqref{Eq0101}
subjects to find a primitive of integral for the left hand side of \eqref{Eq0206}.

Let $\alpha_2=\alpha_3=\alpha_4=\alpha_5=0$,
we get a linear problem $-u''(x)=\alpha_1 u$, $x\in(0,1)$ with $u(0)=u(1)=0$. After some elementary calculation,
one may have the solution
\[u(x) = C_1 \sin\left(\sqrt{\alpha_1}x\right),~~~\alpha_1 = n^2\pi^2,\]
where $n=1,2,3,\cdots$ and $C$ is an arbitrary constant. From the viewpoint of dynamical system,
we transform the second-order ODE into a first-order system by $v(x)=u'(x)$ introduced above. Hence,
\[\frac{du}{dx} = \pm\sqrt{\eta^2-\alpha_1u^2},\]
which admits the solution given by
\[\pm{x} = \frac{1}{\sqrt{\alpha_1}} \left[ \arcsin\left(\frac{\sqrt{\alpha_1} u(x)}{\eta}\right)
- \arcsin\left(\frac{\sqrt{\alpha_1} u(0)}{\eta}\right) \right].\]
Since $u(0)=0$, we can simplify the above expression as $u(x) = \frac{\eta}{\sqrt{\alpha_1}}\sin\left(\sqrt{\alpha_1}x\right)$, 
where we have applied the condition $v(0)=\eta$ to remove the negative sign in $\pm{x}$.
This also implies that for every fixed $\alpha_1=n^2\pi^2>0$, model \eqref{Eq0203} has unique solution
up to a constant factor. However, it shall be an obvious truth that uniqueness of solution will be
violated once the source function on the right hand side becomes nonlinear and satisfies some
mild conditions.

Consider $\alpha_1=\alpha_3=\alpha_4=\alpha_5=0$ in model \eqref{Eq0203}, i.e., $-u''(x)=\alpha_2u^2$.
Then by \eqref{Eq0206} one will have
\begin{equation}\label{Eq0207}
  \int^{u(x)}_{u(0)} \frac{du}{\pm\sqrt{\eta^2 - \frac{2\alpha_2}{3} u^3}} = x.
\end{equation}

Because it has the cubic nonlinear term in the square root, we do not expect an explicit solution
to such integral. With the help of hypergeometric function or elliptic integral, one may get
\begin{equation}\label{Eq0208}
\int^{u(x)}_{u(0)} \frac{du}{\sqrt{\eta^2 - \frac{2\alpha_2}{3} u^3}} =
\frac{u}{\eta}  {}_2F_1\left(\frac{1}{3}, \frac{1}{2}; \frac{4}{3};
\frac{2\alpha_2}{3\eta^2} u^3\right) - \frac{u(0)}{\eta}  {}_2F_1\left(\frac{1}{3}, \frac{1}{2};
\frac{4}{3}; \frac{2\alpha_2}{3\eta^2} u(0)^3\right).
\end{equation}
or
\begin{equation}\label{Eq0209}
\int^{u(x)}_{u(0)} \frac{du}{\sqrt{\eta^2 - \frac{2\alpha_2}{3} u^3}}
= \frac{2}{\sqrt{3}} \left( \frac{3}{2\alpha_2 \eta} \right)^{1/3}
\left[ F\left( \arccos\left( \frac{1+\sqrt{3} \beta x}{1-\sqrt{3} \beta x} \right),
k \right) \right]_{x=u(0)}^{x=u},
\end{equation}
where $F(\varphi, k)$ is the incomplete elliptic integral of the first kind,
$\beta = \left( \frac{2\alpha_2}{3\eta^2} \right)^{1/3}$ and $k^2 = \frac{2+\sqrt{3}}{4}$.

If $\alpha_1=1$, $\alpha_2=\alpha_4=\alpha_5=0$ and $\alpha_3=-1$ in \eqref{Eq0203}, we have
\begin{equation}\label{Eq0210}
  \int^{u(x)}_{u(0)} \frac{du}{\pm\sqrt{\eta^2 - \left(u^2-\frac{1}{2}u^4\right)}} = x.
\end{equation}

Again, with the help of hypergeometric function or elliptic integral, we shall obtain
\begin{equation}\label{Eq0211}
\int^{u(x)}_{u(0)} \frac{du}{\pm\sqrt{\eta^2 - \left(u^2-\frac{1}{2}u^4\right)}}
= \frac{u(x)}{\eta}  {}_2F_1\left( \frac{1}{4}, \frac{1}{2}; \frac{5}{4};
\frac{2u(x)^2 - u(x)^4}{2\eta^2} \right) - \frac{u(0)}{\eta} 
{}_2F_1\left( \frac{1}{4}, \frac{1}{2}; \frac{5}{4}; \frac{2u(0)^2 - u(0)^4}{2\eta^2} \right).
\end{equation}

\begin{equation}\label{Eq0212}
\int^{u(x)}_{u(0)} \frac{du}{\pm\sqrt{\eta^2 - \left(u^2-\frac{1}{2}u^4\right)}}
= \frac{\sqrt{2}}{a} \left[ F\left( \arcsin\left( \frac{a}{u(0)} \right), \frac{b}{a} \right)
- F\left( \arcsin\left( \frac{a}{u(x)} \right), \frac{b}{a} \right) \right].
\end{equation}
where $0<\eta^2 < \frac{1}{2}$, $a=\sqrt{1-\sqrt{1-2\eta^2}}$,
$b=\sqrt{1+\sqrt{1-2\eta^2}}$. For the case $\eta^2 > \frac{1}{2}$,
it may need more tools of complex analysis and hence we choose to omit here. Anyhow, we can see the limit of analysis to identify all solutions. 

\section{Numerical methods}
In this section, we construct a numerical method to identify multiple solutions of the considered problems, which includes the discretization method, how to choose the initial guesses, and how to solve the nonlinear problems.

\subsection{Randomized neural networks}
Neural networks are one important technique to solve forward and inverse partial differential equations. Here we use the randomized neural networks as the tool to discretize the considered problems due to the high efficiency and accuracy \cite{ying2024an}, compared to the fully trained neural networks, and the mesh-free and flexibility in handling the complex domain, compared to the classic numerical methods. For completeness, we give it a brief introduciton.

For a given fully connected neural network, say a shallow neural network $\Theta: \mathbb{R}^{d_{in}} \to  \mathbb{R}^{d_{out}}$ with a width $m$, the network is capable of representing all functions defined by $\Theta(x) = G_2\circ G_{1}\circ G_0(x)$, where $G_0(x)=x$, $G_1=\sigma(w_1 G_{0} +b_1) $, and $G_{2}=w_2G_{1}+b_2$ is the output. Here $\sigma$ is the activation function, which is taken to be $\tanh(x+0.5)-\tanh(x-0.5)$ for the later analysis purpose, $w_1\in \mathbb{R}^{m\times d_{in}}$ and $b_1\in\mathbb{R}^{m}$ are the weights and the bias in the hidden layer, while $w_2\in \mathbb{R}^{d_{out}\times m}$ and $b_2\in\mathbb{R}^{d_{out}}$ are the ones in the output layer. For the randomized neural network, instead of determining all their weights and bias according to the minimization of a given loss function, it randomly determines the weights and bias in the hidden layer (e.g., determine $w_1$ and $b_1$ according to a normal distribution) and only leaves parameters in the last layer tunable. Without loss of generality, we assume the bias term $b_2$ in the last layer to be zero.

Then in order to find a solution of equation \eqref{Eq0101} based on neural networks, we construct the following loss function and the minimization problem
\begin{equation}
\label{lsp}
\min_{ u \in H^{2}({\Omega})} L( u)  =\min_{ u \in H^{2}({\Omega})} \sum_{j=1}^2 \beta^2_j L_j,
\end{equation}
where the positive coefficients $\beta_{j}$ serve as penalty weights, and the residual terms $L_{j}$ measure the errors within the domain and on the boundary respectively, which are defined by
\[
L_{1}({u}) = \int_{\Omega}| - \Delta {u} (x)-\lambda f(x,u(x))|^{2}dx,\quad L_{2}({u}) = \int_{\partial \Omega}\bigl|{u} (\mathbf{s}) - g(\mathbf{s})\bigr|^{2}d\mathbf{s}.
\]

By approximating the function ${u}(x)$ with a shallow network of width $m$ and employing Monte Carlo sampling to discretize the loss, we then have the following discrete minimization problem
\begin{equation}
\label{lspm_d}
{u_{md}^*} = \arg\min_{ u_m \in \mathcal{H}_m}  L^{d}( u_m) =\arg\min_{ u_m \in \mathcal{H}_m} \sum_{j=1}^2 \beta^2_j l_j,
\end{equation}
where $\mathcal{H}_m$ denotes the hypothesis space, $l_1$ and $l_2$ are the discrete forms given by
\begin{equation}
\nonumber
    \begin{aligned}
        l_1( u_m) &= \frac{|\Omega|}{m_{\Omega}}\sum_{j=1}^{m_{\Omega}} \left| - \Delta {u_m} (x_j)-\lambda f(x_j,u_m(x_j)) \right|^2,   x_j\in\Omega, \\
        l_2( u_m)  &=\frac{|\partial\Omega|}{m_{\partial\Omega}} \sum_{j=1}^{m_{\partial\Omega}} \left|    u_m \left( \s_j\right) - g(\s_j)  \right|^2, \quad \bold{s}_j\in\partial\Omega
    \end{aligned}
\end{equation}
with $m_{\Omega}$ and $m_{\partial\Omega}$ being the numbers of (uniformly distributed) sampling points in $\Omega$ and on $\partial\Omega$, respectively.

For a fully connected network, what we need to do is invoking an optimization method  to minimize the loss function to obtain a solution. However, for a randomized neural network, due to the special structure of the hypothesis space $\mathcal{H}_m=\{\phi_1, \cdots, \phi_m\}$, it is equivalent to the following nonlinear least square problem:
\begin{equation}
\label{min_equiv_nonlinear}
\min_{  \tilde{u}\in\mathcal{H}_m  }\left\|N(\tilde{u}) - b\right\|^2,
\end{equation}
where $\tilde{u}=\sum_{k=1}^m\tilde{u}_k\phi_k \in \mathcal{H}_m$,
\begin{equation}\nonumber
N(\tilde{u}) =  \begin{matrix} 
  \left [ \begin{matrix} 
    \beta_{1} \left[ -\Delta  \tilde{u}(x_j) - \lambda f(x_j, \tilde{u}(x_j)) \right]_{m_{\Omega}\times 1} \\
       \beta_{2}\left[   \tilde{u}(x_j)   \right]_{m_{\partial\Omega}\times 1 }\\
   \end{matrix} \right]   \\
   \end{matrix},
\end{equation}
and
\begin{equation}\nonumber
b =  \begin{matrix} 
  \left [ \begin{matrix} 
    \beta_{1} \left[ 0 \right]_{m_{\Omega}\times 1} \\
       \beta_{2}\left[   g(\s_j)   \right]_{m_{\partial\Omega}\times 1 }\\
   \end{matrix} \right]   \\
   \end{matrix}.
\end{equation}
To solve this nonlinear least square problem, we can call some nonlinear solver, like the function \texttt{least}$\_$\texttt{square} from the \texttt{scipy} package, by solving the normal equations iteratively to get an approximation. Of course, the boundary part in \eqref{min_equiv_nonlinear} can be ignored if we use the hard constraint instead of the soft constraint \cite{lai2025hard}. Therefore, with a little abuse of notations and for easy presentations in the following, without loss of generality, we will assume the boundary condition is exactly satisfied by the neural network.
Next, we are going to construct new method to compute multiple solutions of semilinear problem from two aspects: initializations and optimizations.

\subsection{Initializations}
For a nonlinear problem, it is critical to have a sufficient good initial guess to guarantee the convergence using some Newton-like methods. To find multiple solutions, it is exactly the same here. Meanwhile, we expect that if the considered problem has more than one solution, each solution should have its own convergence region so that different initializations would lead to different final solutions. As shown in \cite{zou2025learning}, random initializations can be used to uncover multiple solutions. So in this subsection, we like to propose a new method to efficiently construct different initial guesses to eventually find different solutions effectively.

Given a randomized neural network with the width $m$, the sampling points $\{x_j\}^{m_{\Omega}}_{j=1}$ uniformly distributed in $\Omega$ (no sampling points on the domain boundary as the Dirichlet boundary condition is exactly satisfied by the assumption), and an arbitrary initial guess $u^{(0)}_{1}$, by solving the nonlinear least square problem \eqref{min_equiv_nonlinear}, suppose we have obtained one solution $u_1$. Next, with this fixed shallow neural network and the fixed sampling points, to obtain different solutions, we propose a new method to systematically design how to choose initial guess properly to fulfill this goal.

First of all, we introduce the following discretized characteristic matrix
\[\Phi=[\Phi_{ij}]_{m_{\Omega}\times m}\mbox{ with }\Phi_{ij}=\phi_j(x_i),\]
which is used to judge how good the coefficient vector of $\zeta\equiv(\tilde{u}_1, \cdots, \tilde{u}_m)$ in the linear combination of $\Phi \zeta$ is to fit the problem \eqref{min_equiv_nonlinear}. Then for the solution $u_1$, we have $\left(u_1(x_1), \cdots, u_1(x_{m_{\Omega}})\right)^T=\Phi \zeta_1$ for some coefficient vector $\zeta_1$. To effectively construct the initial guess, we will follow the idea of the Newton method and firstly consider the nonlinear problem around the (first) found solution $\Phi \zeta_1$ locally.

To do that, we assume the following discretized nonlinear residual function
\begin{equation}\label{F-equa}
F(\zeta)\equiv-\Delta \Phi\zeta-\lambda f(x,\Phi\zeta)
\end{equation}
 with $f(x,\Phi\zeta)\equiv\left( f(x_j, (\Phi\zeta)_j )\right)_{m_{\Omega\times 1}}$ and $\Delta \Phi\equiv \left( \Delta \phi_j(x_i) \right)_{m_{\Omega}\times m}$. Then we know that $F(\zeta_1)=0$ and
\[F(\zeta_1+\delta\zeta)\approx F(\zeta_1)+J(\zeta_1)\delta\zeta =J(\zeta_1)\delta\zeta\]
with $J(\zeta_1)$ denoting the Jacobian matrix at the point $\zeta_1$ and $\delta\zeta$ meaning a small perturbation. This indicates the local perturbation of the residual $\delta F$ is approximately determined by the quantity $J(\zeta_1)\delta\zeta$. Therefore, if there is another solution near $\zeta_1$, the search direction $\delta\zeta$ with the norm of $J(\zeta_1)\delta\zeta$ being relatively small would be a good choice, since small residual in some sense indicates being sufficient close to a solution of the considered problem. For the notation simplicity, we call this kind of search direction to be weakly variant search direction (WVSD). As for WVSD, we have the following easy result to tell us how to choose these directions.
\begin{thm}
For $\delta F\equiv J(\zeta_1)\delta\zeta$, the best WVSD is any right eigenvector corresponding to the smallest singular value.
\end{thm}
\begin{proof}
For the matrix $J(\zeta_1)$, we have the singular value decomposition (SVD) as follows
\[J(\zeta_1)\delta\zeta=U_J\Sigma_J V^T_J,\]
where $U_J$ and $V_J$ are the orthonormal matrices containing the left and right eigenvectors, and $\Sigma_J$ is the diagonal matrix with singular values $\lambda_1\geq \lambda_2\geq \cdots \lambda_m\geq 0$.

Let $\delta\zeta=V_JK$ for some vector $K=(k_1,\cdots,k_m)$, and we have
\[\|\delta F\|=\|J(\zeta_1)\delta\zeta\|=\|U_J\Sigma_J K\|=\|\Sigma_J K\|= \sqrt{\sum_{j=1}^m\lambda^2_jk^2_j},\]
from which we know that $\min_{\|\delta\zeta\|=1} \|\delta F\|= \lambda_m,$
exactly when we take $\delta\zeta$ to be right eigenvector in $V_j$ corresponding to the singular value $\lambda_m$. Therefore, the conclusion follows.
\end{proof}

Based on the solution representation form of $\Phi \zeta$, it means finding multiple solutions is just equal to uncover different coefficient vectors $\zeta$ for the nonlinear least square problem \eqref{min_equiv_nonlinear}. However, due to the randomness in generating these basis functions $\phi_i$, it may lead to the column vectors of the characteristic matrix $\Phi$ is (weakly) linear dependent so that even for some perturbation $\delta\zeta$ with large norm, it is possible to have
$$\Phi\zeta\approx\Phi\left(\zeta+\delta\zeta \right)=\Phi\zeta+\Phi\delta\zeta,$$
which is the easy consequence when $\Sigma V\delta\zeta\approx 0$ with the SVD of $\Phi$ being $U\Sigma V$ and the smallest singular values close to zero (i.e., $\delta\zeta$ is taken to be the right eigenvector corresponding to the smallest singular value). Therefore, it is proper to introduce the following optimization subproblem:
\[
\left\{
\begin{aligned}
& \min_{\mathbb{R}^m} \|J(\zeta_1)\delta\zeta\| \\
& \text{s.t.} \quad \|\Phi\delta\zeta\| = 1
\end{aligned}
\right.
\]
to determine the possible (local) search direction for the purpose of identifying effective initial guesses.

Based on the SVD, it turns out it is unnecessary to solve this subproblem. However, we need to do some modifications in order to improve the numerical stability. First of all, for the characteristic matrix $\Phi=U\Sigma V$, since the search direction $\delta\zeta$ is not changed by a multiplication of a constant, the constraint of $\|\Phi\delta\zeta\| = 1$ can be imposed or fulfilled by the normalization technique (i.e., $\delta\zeta/\|\Phi\delta\zeta\| $). Hence in order to avoid the overflow problems caused by the extremely small values of $\|\Phi\delta\zeta\|$, we truncate the SVD decomposition so that the singular values less than or equal to a given threshold $\epsilon_1>0$ are neglected. By the truncation, we then have
$\Phi_t=U_t\Sigma_t V^T_t,$
where $U_t\in\mathbb{R}^{m_{\Omega}\times r}$, $\Sigma_t\in\mathbb{R}^{r\times r}$, and $V_t\in\mathbb{R}^{m\times r}$ are the matrices consisting of the left eigenvectors, the singular values larger than $\epsilon_1$, and the right eigenvectors correspondingly.  Meanwhile, in order to minimize the effect of different numbers of sampling points, following the definition of the root mean square, the norm to measure the variation in terms of the solution is redefined as follows
\begin{equation}
\label{h_norm}
\|\delta\zeta\|_h\equiv\sqrt{\frac{1}{m_{\Omega}}\|\Phi_t\delta\zeta\|^2 },
\end{equation}
With this change, the constraint $\|\Phi\delta\zeta\|=1$ in the minimization problem will be changed to $\|\delta\zeta\|_h=1$ correspondingly. Furthermore, to reduce the computational costs, we can also setup another threshold $\epsilon_2>0$ such that the potential direction is only searched in a relatively small subset spanned by the right eigenvectors of small eigenvalues of the Jacobian matrix $J(\zeta_1)$. This is also essential for increasing the robustness of the method. As for the minimization problem, the exact solution is the right eigenvector of the smallest singular value multiplied by certain constant to satisfy the equality constraint. Indeed, this is the best choice of the search direction locally. However, from the global point of view, it could be better to include also the right eigenvectors corresponding to the small eigenvalues (e.g., less than the given $\epsilon_2$) as the choices of the search direction. For the latter case with multiple eigenvectors, we actually just need to set each eigenvector as the candidate of the search direction $\delta\zeta$ since what we are interested in the search direction for constructing effective initial guesses.

Noticing that $\|\Phi_t\delta\zeta\|=\|\Sigma_tV^T_t\delta\zeta\|$ as the columns of $U$ are orthonormal with each other, and the singular values $\sigma_1, \cdots, \sigma_r$ in $\Sigma_t$ after the truncation are positive, if we take $\delta\zeta\equiv Wz$ for some vector $z\in\mathbb{R}^r$ with
$W\equiv \sqrt{m_{\Omega}} V_t \mbox{diag}(\frac{1}{\sigma_1},\cdots,\frac{1}{\sigma_r}),$
it implies that
$$\|\delta\zeta\|_h=\|Wz\|_h=\sqrt{\frac{1}{m_{\Omega}}\|\Phi_tWz\|^2}=\sqrt{\|\Sigma_tV_t^T V_t\mbox{diag}(\frac{1}{\sigma_1},\cdots,\frac{1}{\sigma_r}z)\|^2}=\|z\|.$$
Therefore, for the constraint $\|\delta\zeta\|_h=1$, we have $\|z\|=1$ when $\delta\zeta=Wz$. So we consider the SVD of the modified Jacobian $J(\eta_1)W\equiv U_J\Sigma_J V^T_J$ and once the $\epsilon_2>0$ is given,  we can find the right eigenvectors $v_J^{(k)}$ from $V_J$ whose singular values are smaller than $\epsilon_2$, followed by the multiplication $Wv_J^{(k)}$ to determine the potential search directions together with arbitrarily chosen step sizes.

For clarity purpose, we summarize the algorithm of the initialization in the following.
%
%
%
	\begin{center}
		\begin{minipage}{0.94\linewidth}
			\hrule
			\vspace{0.5em}
			\noindent\textbf{Algorithm 1: Initializations }
			\vspace{0.5em}
			
			\noindent\textbf{Input:}
			A randomized neural network with $m$ hidden layers, $m_{\Omega}$ sampling points, and a precomputed solution $u_1$ with the coefficient vector $\zeta_1$. Meanwhile, specify $\epsilon_1$ and $\epsilon_2$.
			
			\begin{enumerate}
				\item For the characteristic matrix $\Phi$ and the given $\epsilon_1$, construct the truncated SVD $\Phi_t=U_t\Sigma_tV^T_t$ and the $W$ matrix.
				\item For the nonlinear problem $F(\zeta)$, compute the Jacobian matrix $J(\zeta_1)$ at the point $\zeta_1$.	
				\item Compute the SVD of the modified Jacobian matrix $J(\zeta_1)W=U_J\Sigma_J V^T_J.$
				\item Given $\epsilon_2$, find the right eigenvectors in $V_J$ corresponding to the singular eigenvalues less than $\epsilon_2$.
				\item For each $v^{(k)}_J$ of the right eigenvectors, we construct the search direction by $Wv^{(k)}_J$.
			\end{enumerate}
			
			\noindent\textbf{Output:}
			the initial guesses constructed by $\zeta=\zeta_1\pm \alpha_{init} Wv^{(k)}_J$ with the step size $ \alpha_{init}$.
			\vspace{0.5em}
			\hrule
		\end{minipage}
	\end{center}		
In practice, we take $\epsilon_1$ to be the largest module of the eigenvalues multiplied by a given small constant. Meanwhile, we can consider the right eigenvectors and thus the search directions corresponding to the first few $K$ smallest singular values without using $\epsilon_2$. And the step size $ \alpha_{init}$ can be chosen empirically. By default, we set $\alpha_{init}=0.5, 1, 2, 4, 8$ and set $K=6$ in the numerical tests.

\subsection{Optimization method}
This part is devoted to the optimization method based on the deflation technique for the discretized problem $F(\zeta)$ given in \eqref{F-equa}. Assume currently we already find some solutions, say $\zeta_1, \cdots, \zeta_{n-1}$ for some integer $n$. With the initial guess $\zeta_n^{(0)}$ produced in Algorithm 1, we try to find a new solution $\zeta_n \notin U_{s}\equiv\{\zeta_1,\cdots, \zeta_{n-1}\}$.

For this purpose, first of all, we define the deflation operator $\mu(x; U_{s})$ so that new solution is expected to find for the considered problem. Commonly used deflation operators for a single root (i.e., $U_s=\{\zeta_1\}$) includes
\begin{enumerate}
\item[$\bullet$]  $\displaystyle\mu(x; \{\zeta_1\})=\frac{1}{\|x-\zeta_1\|}$;\\
\item[$\bullet$] $\displaystyle \mu(x; \{\zeta_1\})=\frac{1}{\|x-\zeta_1\|^{\theta}}+\sigma$ with $\theta=2$ and $\sigma=1$;
\item[$\bullet$] $\displaystyle \mu(x; \{\zeta_1\})=\exp\left(\frac{1}{\|x-\zeta_1\|}\right)$.
\end{enumerate}
Here the second deflation operator includes the $\theta$ parameter to overcome the higher multiplicity of a root, while the last exponential one is proposed to incorporate all possible values of $\theta$ instead of choosing the parameter value empirically. Then extended from the single root to the multiple roots, we have the following deflation operator
\begin{equation}
\label{df-m}
 \mu(x; U_{s})=\Pi_{\zeta_k\in U_s} \mu(x; \zeta_k).
\end{equation}

Then considering the following least square problem
\begin{equation}
\label{min-p}
\min_{\zeta\in\mathbb{R}^m}\frac{1}{2}\|F(\zeta)\|^2,
\end{equation}
we adapt the good deflated Gauss-Newton method proposed in \cite{riley2025deflation}, where the norm $\|\cdot\|_h$ is used in the deflation to measure the effective distance of differnet solutions. Its novelties and advantages lie in the following three aspects:
\begin{enumerate}
\item[$\bullet$] The exponential deflation operator is used;
\item[$\bullet$] Instead of solving the nonlinear linear square problem, like the Levenberg–Marquardt method, only the linearized part is considered to design the deflation algorithm, which avoids the calculation of the Hessian matrix and thus is computationally cheaper. Meanwhile, it can avoid the convergence to the stationary points;
\item[$\bullet$] Theoretically it guarantees (see \cite[Theorem 3.5]{riley2025deflation}) that the local convergence of the proposed method (i.e., converge to the local minima that are of interest), and it would not converge to any of the deflated points (i.e., the founded solution in $U_s$) under some requirements of the deflation operator and the Jacobian matrix.
\end{enumerate}

For clarity and completeness purpose, we present the good deflation Gauss-Newton algorithm in the following.
	\begin{center}
		\begin{minipage}{0.98\linewidth}
			\hrule
			\vspace{0.5em}
			\noindent\textbf{Algorithm 2: Good deflated Gauss-Newton algorithm (GDGN)}
			\vspace{0.5em}
			
			\noindent\textbf{Input:}
			Given an initial guess $\zeta_n^{(0)}$, and the set of founded solutions $U_s$ as well as the corresponding deflation operator $\mu(\zeta; U_s)$ defined in \eqref{df-m}, where the norm $\|\cdot\|_h$ defined in \eqref{h_norm} is used. Meanwhile, choose a $\epsilon\in[0,1]$, a $tol$ and define $\eta(\zeta)=\log(\mu(\zeta; U_s))$.
			
				\begin{enumerate}
				\item[] \textbf{While} $\|p^{(k)}\|>tol$ \textbf{do}
				\begin{enumerate}
				\item[] $p^{(k)}=\arg\min_p \|F(\zeta^{(k)}_n)+J_F(\zeta^{(k)}_n)p\|$ with $J_F(\zeta)$ being the Jacobian of the minimized function given \eqref{min-p} in at the point $\zeta$.
				\item[] \textbf{if} $\langle p^{(k)}, \nabla\eta(\zeta_n^{(k)}) \rangle>\epsilon$ \textbf{then}
				\begin{enumerate}
				\item[] $\beta=1 - \langle p^{(k)}, \nabla\eta(\zeta_n^{(k)}) \rangle$;
				\item[] $\zeta_n^{(k+1)} = \zeta_n^{(k)} + \beta^{-1}p^{(k)}$;
				\end{enumerate}
				\item[] \textbf{else}
				\begin{enumerate}
				\item[] $\zeta_n^{(k+1)} = \zeta_n^{(k)} + \alpha p^{(k)}$ with $\alpha$ being determined by the quadratic line search on the nonlinear problem defined in \eqref{min-p};
				\end{enumerate}
			\item[] \textbf{end if}
			\end{enumerate}
			\item[] \textbf{end while}
			\end{enumerate}
			
			\noindent\textbf{Output:}
			if finding a new solution $\zeta_n$, update $U_s=U_s\cup \{\zeta_n\}$.
			\vspace{0.5em}
			\hrule
		\end{minipage}
	\end{center}	
In this algorithm, for the sub-minimization problem $\min_{p\in\mathbb{R}^m} \|F(\zeta^{(k)}_n)+J_F(\zeta^{(k)}_n)p\|,$ since the sampling points $m_{\Omega}$ is generally larger than the width of the hidden layer $m$, so the minimizer can be obtained by solving the following
\[J^T_F(\zeta^{(k)}_n)J_F(\zeta^{(k)}_n) p = - J^T_F(\zeta^{(k)}_n)F(\zeta^{(k)}_n),\]
which is efficiently solved based on the QR decomposition or the SVD decomposition of $J_F(\zeta^{(k)}_n)$. Meanwhile, in practice, to see a $\zeta_n$ is a solution or not, we setup another threshold $R_{\epsilon}>0$ and we accept the approximation $\zeta_n$ as a solution if the norm of the residue is less than or equal to $R_{\epsilon}$, where $R_{\epsilon}$ will be specified in the numerical tests.

\subsection{Local convergence analysis}
Here we present the local convergence of the used method for solving the semilinear problem \eqref{Eq0101} for some particular solution $u^*\in H^2(\Omega)$. Without loss of generality, we assume $g=0$ and $\Omega\subset\mathbb R^n$ is a bounded $C^{1,1}$ domain with $n\leq 3$ for simplicity. Meanwhile, we also assume the following assumptions:
\begin{enumerate}
\item[A.] The function $f(x,s)$ is measurable in $x$ and
twice continuously differentiable in $s$ for
almost every $x\in\Omega$.
Meanwhile, there exists a constant $M>0$ such that,
in a neighborhood $B_{u^*,\delta}$ of $u^*$, $
|f_u(x,s)|+|f_{uu}(x,s)|\leq M
$
for almost every $x\in\Omega$ and $s\in B_{u^*,\delta}$.
\item[B.] [Non-degeneracy]
The linearized homogeneous Dirichlet problem
\[
\begin{cases}
-\Delta h-\lambda f_u(x,u^*)h=0,
&x\in\Omega,\\
h=0,&x\in\partial\Omega
\end{cases}
\]
has only the trivial solution. Equivalently, $
\ker L_*=\{0\}$ with $L_*h=F'(v^*)h
=-\Delta h-\lambda f_u(x,u^*)h.$
\end{enumerate}
\begin{remark}
The non-degeneracy condition is satisfied by some solutions, like the positive solutions of the two-dimensional Lane-Emden problem  \cite{grossi2022non}, which we will investigate in the numerical section.
\end{remark}

Based on these two assumptions, we firstly have the following.
\begin{thm}\label{linear_invertible}
Suppose $X=H^2(\Omega)\cap H_0^1(\Omega)$ and $
Y=L^2(\Omega)$, then we have $L_*:X\to Y$ is an isomorphism and
\[
\|h\|_{H^2(\Omega)}
\lesssim \|L_*h\|_{L^2(\Omega)}
\qquad\forall h\in X. \]
\end{thm}
\begin{proof}
For the Laplace operator $A=-\Delta$ defined on $X$, by the Lax--Milgram theorem and the elliptic
regularity $C^{1,1}$ on $\Omega$, it is easy to see that $
A: X\to Y$
is bijective and therefore we know $
A^{-1}
$ exists and
is bounded.

Let $
V(x)=\lambda f_u(x,u^*)$. By Assumption A, we have $
V\in L^\infty(\Omega).
$
Denoting by $M_V$ the multiplication operator
\[
M_Vh=Vh,
\]
then we obtain $L_*=A-M_V=A(I-A^{-1}M_V)$. Since $A^{-1}$ is bounded and the embedding $
H^2(\Omega)\hookrightarrow L^2(\Omega)
$
is compact by the Rellich–Kondrachov theorem, $A^{-1}: L^2(\Omega) \to L^2(\Omega)$ is a compact operator. Consequently, $K=A^{-1}M_V:L^2(\Omega)\to L^2(\Omega)$
is also compact as $M_V$ is linear bounded.

Therefore, by the Riesz--Schauder theorem and the assumption B, $I-A^{-1}M_V$ is
invertible on $L^2(\Omega)$. Since $A$ is
invertible, we hence know that $L_*:X\to Y$ is an isomorphism, from which we obtain the inequality.
\end{proof}

For the nonlinear problem \eqref{Eq0101}, corresponding we can define the nonlinear operator $F:X\to Y$, where
\begin{equation}\label{eq:F}
F(v)=-\Delta v-\lambda f(x,v).
\end{equation}
Then for $\forall h\in X$ satisfying $u^*+th \in B_{u^*,\delta}$ for $t\in[0,1]$ and by the Taylor expansion, we have $F(t)\equiv F(u^*+th)=F(0)+F'(0)t + \int_{0}^t F^{''}(T)(t-T)dT$, from which we have
 \[ F(u^*+h)= F(1)=L^*h + R(h)\]
with $R(h)=-\lambda\int_0^1(1-t)f_{uu}(x,u^*+th)h^2 dt.$ By the assumption A and the Sobolev embedding $H^2(\Omega)\hookrightarrow C^{0,\alpha}(\Omega) \hookrightarrow L^{\infty}(\Omega)$ for $n\leq 3$, we have
\begin{align*}
\|R(h)\|^2_{L^2}=\frac{\lambda^2 M^2}{2}\int_{\Omega}h^4dx\lesssim\|h\|^2_{L^{\infty}}\int_{\Omega}h^2dx=\|h\|^2_{L^{\infty}} \|h\|_{L^2}^2\lesssim  \|h\|^4_{H^2},
\end{align*}
which give the local estimate of the nonlinear remainder term. Based on the analysis, then we are able to present the following local stability result as a easy consequence.
\begin{thm}
For $\forall u\in X$ satisfying $\|u-u^*\|_{H^2}<r$ (i.e., $u \in B_{u^*,r}$) for some $0<r\leq \delta$,
we have
\begin{equation}\label{eq:local-stability}
\|u-u^*\|_{H^2}
\lesssim \|F(u)\|_{L^2}.
\end{equation}
\end{thm}

\begin{proof}
Set $h=u-u^*$, and we then have
\[
F(u)=F(u^*+h)=L_*h+R(h).
\]
For $L_*h=F(u)-R(h)$ with $h\in B_{u^*,\delta}$, by Theorem \ref{linear_invertible} and the previous calculation, we have
\[
\begin{aligned}
\|h\|_{H^2} \lesssim \|L_*h\|_{L^2} \leq \|F(u)\|_{L^2}+\|R(h)\|_{L^2}\lesssim \|F(u)\|_{L^2}+\|h\|^2_{H^2}\equiv C\left( \|F(u)\|_{L^2}+\|h\|^2_{H^2}\right),
\end{aligned}
\]
where the generic constant is denoted by $C$.

We choose $0<r\leq \delta$ such that $Cr\leq\frac12$. Then if $\|h\|_{H^2}<r$, we have $C\|h\|_{H^2}^2
\leq\frac12\|h\|_{H^2}$
and $
\frac12\|h\|_{H^2}
\leq C\|F(u)\|_{L^2}
$. Therefore,
\[
\|u-u^*\|_{H^2}
\leq 2C\|F(u)\|_{L^2},
\]
which concludes the theorem with the generic constant $2C$ in the inequality \eqref{eq:local-stability}.
\end{proof}

Based on the stability result and following the routine analysis we commonly used \cite{ying2024an,li2025neural,li2026shallow}, then we are able to prove the error  estimate without considering optimization errors in solving the nonlinear least square problem. Assume the boundary condition is satisfied exactly, then $\|F(u)\|_{L^2}=L^1(u)=L(u)$ defined in \eqref{lsp} when $\beta_1=1$, and similarly the corresponding discretized problem is the minimization problem $L^d$ given in \eqref{lspm_d} without the boundary term neither. Then, $\forall u\in \mathcal{H}_m$,
\begin{align*}
L(u^*_{md})&=L(u^*_{md})-L(u^*)\leq L(u^*_{md}) - L^d(u^*_{md}) + L^d(u) - L(u) + L(u) - L(u^*)\\
&\leq \inf_{u\in\mathcal{H}_m}\left|L(u) - L(u^*)\right| + 2\sup_{u\in \mathcal{H}_m} \left| L(u) - L^d(u) \right|.
\end{align*}
Meanwhile, noticing that $F(u^*)=0$, we have
\begin{align*}
L(u) - L(u^*) &= \|F(u)\|_{L^2}^2\leq \left( \|F(u)-F(u^*)\|_{L^2}+\|F(u^*)\|_{L^2}\right)^2\leq \|F(u)-F(u^*)\|_{L^2}^2 \\
&=\int_{\Omega} \left( -\Delta u-\lambda f(x,u) +\Delta u^*+\lambda f(x,u^*) \right)^2dx\\
&\leq 2 \int_{\Omega} \left( -\Delta u +\Delta u^* \right)^2dx + 2\int_{\Omega} \left( \lambda f(x,u) -\lambda f(x,u^*) \right)^2dx\\
&\lesssim \|u-u^*\|^2_{H^2} + \|u-u^*\|^2_{L^2} \leq \|u-u^*\|^2_{H^2},
\end{align*}
where the last second inequality is obtained due to the Lipschitz continuity of the nonlinear equation $f$ with respect to $u$ as assumed in the assumption A. Therefore, we have
\[
L(u^*_{md})\lesssim \inf_{u\in\mathcal{H}_m}\|u-u^*\|^2_{H^2} + \sup_{u\in \mathcal{H}_m} \left| L(u) - L^d(u) \right|,
\]
where the first term characterize the approximation error, and the second term is commonly known as the statistical or generalization error.
As for the approximation error, we have the following result proved in \cite[Corollary 3.8]{bi2026adaptive}.
\begin{thm}
Assume that $u^*$ have a global extension $u^*_e$ such that $u^*_e\in L^1(\mathbb{R}^n)$ and its Fourier transform $\widehat{u^*_e}\in L^1(\mathbb{R}^n)$. For the randomized neural networks with random variables $(w_1,b_1)$ sampled independent and identically distributed from the uniform distributions over the bounded domain $\{ \xi\in \mathbb{R}^n~|~\|\xi\|\leq \frac{M}{2R}\}\times[-M, M]$, if the parameter $M$ is chosen as
$M=O(m^{\frac{1}{n+1+2\eta}})$ and $u^*\in H^{2+k}\subset H^2$, then we have
\[\mathbb{E}\left[ u^* - u^{w_1,b_1}_{w_2} \right] \lesssim \|u^*\|_{H^4} m^{-\frac{\eta}{n+2\eta +1}},\]
where $k> \frac{n}{2}$ is a positive integer, $\eta=k-\displaystyle\frac{n}{2}$, $R=\max_{x\in\Omega}|x|$,  $u^{w_1,b_1}_{w_2}=\sum_{j=1}^m w_{2,j}\sigma(w_{1,j}x+b_{1,j})$ with $w_{1,j}$ being the $j$-th row of the matrix, and $w_{2,j}$ and $b_{1,j}$ denoting the $j$-th component of the vectors, respectively.
\end{thm}
\begin{remark}
Fix the way to randomly generate the parameters in the hidden layer, by considering the almost optimal way in the analysis to determine the weights in the output layer, we try to show the smallest bound for all possible different choices of parameters in the hidden layer,  from which we know it is reasonable assume the output weights have bounded Euclidean norm under some additional regularity assumption, a commonly-used assumption in the establishing the statistical error.

To do that, suppose $\mathcal{H}_m=\{\phi_1, \cdots, \phi_m\}$ with $\phi_i=w_{2,i}\sigma(w_{1,i}x+b_{1,i})$, where $w_{1,i}\in\{ \xi\in \mathbb{R}^n~|~\|\xi\|\leq \frac{M}{2R}\}$ and $b_{1,i}\in[-M, M]$.
As shown in \cite{bi2026adaptive}, we know
\[w_2 = (w_{2,1},\cdots, w_{2,m})^T, \quad w_{2,k}= \frac{1}{m} \frac{\beta(w_{1,k},b_{1,k})}{\pi(w_{1,k},b_{1,k})} \]
where the function $\beta: \mathbb{R}^n\times \mathbb{R}\to\mathbb{R}$ is defined to satisfy integral representation of the shallow neural network, and $\pi$ is the density function of the uniform distribution satisfying $\frac{1}{\pi}\lesssim M^{n+1}$.

Since $\beta(\xi,s)=\frac{a^n}{(2\pi)^n}\widehat{\sigma}(a) 
\widehat{{u_e^*}}(a\xi)e^{-ias}$ for some constant $a$, we know that it is bounded as ${u^*_e}\in L^1(\mathbb{R}^n)$. Therefore, we have
\begin{align*}
\|w_2\|=\sqrt{\sum_{k=1}^m |w_{2,k}|^2}\lesssim \frac{1}{\sqrt{m}}M^{n+1}
\lesssim m^{-\frac{1}{2}+\frac{n+1}{n+2\eta +1}}.
\end{align*}
This implies that for solution with high regularity, $\|w_2\|\lesssim m^{\alpha}$ for some $\alpha<0$, and increasing number of neurons leads to the decay of the almost optimal output weights. So assuming the solution $u^*\in H^5$ at least or higher, we can reasonably assume the output weights have bounded Euclidean norm independent of the width $m$ for simplicity, which is also required by the local error estimate.
\end{remark}

As for the statistical error, after the set $\mathcal{H}_m$ of basis functions are determined by randomly chosen according to the uniform distributions, we can prove the following.
\begin{thm}
For a randomized shallow network with bounded output weights in the Euclidean norm, we have
\[ \mathbb{E}_{x_1,\cdots,x_{m_{\Omega}}}\left[ \sup_{{u}\in\mathcal{H}_m} \left| L^{d}({u}) - L({u}) \right | \right] \leq C_mm^{-\frac{1}{2}}_{\Omega} \]
for some constant $C_m=O(m^{\frac{3}{2}+\frac{3}{n+1+2\eta}})$.
\end{thm}
\begin{proof}
Here we mainly use the Rademacher complexity to derive the result. Firstly,
\begin{align*}
\mathbb{E}_{x_1,\cdots,x_{m_{\Omega}}}\left[ \sup_{{u}\in\mathcal{H}_m} \left| L^{d}({u}) - L({u}) \right | \right] \leq 2R_{m_{\Omega}}\left(\psi\circ \left(\Delta\mathcal{H}_m+\lambda f\circ\mathcal{H}_m\right)\right),
\end{align*}
where $R_N(\mathcal{F})= \mathbb{E}_{\epsilon_1,\cdots,\epsilon_{N}}\mathbb{E}_{x_1,\cdots,x_{N}}\left[  \sup_{h\in\mathcal{F}}\frac{1}{N} \sum_{i=1}^N\epsilon_ih(x_i) \right]$
with $\epsilon_1,\cdots,\epsilon_{N}$ being the Rademacher random variables, $f\circ\mathcal{H}_m=\{f\circ u(x)\equiv f(x, u(x)) ~|~u\in\mathcal{H}_m \}$, and $\psi(x) = x^2$.

For a basis function $\phi_j(x)=w_{2,j}\sigma(w_{1,j}x+b_{1,j}) \in \mathcal{H}_m$, $\Delta \phi_j=w_{2,j}\|w_{1,j}\|^2\sigma^{''}(w_{1,j}x+b_{1,j})$, from which we know $\sup_{v\in \Delta\mathcal{H}_m}\|v\|_{L^{\infty}}\lesssim mM^{2}$.
Since $\lambda$ is a fixed constant, $\sigma$ is the difference of two $\tanh$ functions, and $f$ is Lipschitz continuous with respect to $u$, we know  $\sup_{v\in \Delta  f\circ \mathcal{H}_m}\|v\|_{L^{\infty}}\lesssim m$. Therefore, we further have
\[R_{m_{\Omega}}\left(\psi\circ \left(\Delta\mathcal{H}_m+\lambda f\circ\mathcal{H}_m\right)\right)\lesssim mM^{2}  R_{m_{\Omega}}\left(\Delta\mathcal{H}_m+\lambda f\circ\mathcal{H}_m\right).\]
Separating by the linearity of  the Rademacher complexity, then for the first term we have
\begin{align*}
R_{m_{\Omega}}\left(\Delta\mathcal{H}_m\right)&\leq\mathbb{E}_{\epsilon_1,\cdots,\epsilon_{m_{\Omega}}}\mathbb{E}_{x_1,\cdots,x_{m_{\Omega}}}\left[  \sup_{\mathcal{H}_m}\frac{1}{m_{\Omega}} \left| \sum_{i=1}^{m_{\Omega}}\epsilon_i\sum_{k=1}^mc_k\Delta\phi_k(x_i)\right| \right]  \leq \frac{1}{m_{\Omega}}\mathbb{E}_{\epsilon_1,\cdots,\epsilon_{m_{\Omega}}}\mathbb{E}_{x_1,\cdots,x_{m_{\Omega}}}\left[ \sqrt{\sum_{k=1}^{m} \left(\sum_{i=1}^{m_{\Omega}}\Delta\epsilon_i\phi_k(x_i)\right)^2 } \right]\\
& \leq \frac{1}{m_{\Omega}}\left(\mathbb{E}_{\epsilon_1,\cdots,\epsilon_{m_{\Omega}}}\mathbb{E}_{x_1,\cdots,x_{m_{\Omega}}}\left[ \sum_{k=1}^{m} \left(\sum_{i=1}^{m_{\Omega}}\epsilon_i\Delta\phi_k(x_i)\right)^2 \right]\right)^{\frac{1}{2}}= \frac{1}{m_{\Omega}}\left(\sum_{k=1}^{m}\mathbb{E}_{\epsilon_1,\cdots,\epsilon_{m_{\Omega}}}\mathbb{E}_{x_1,\cdots,x_{m_{\Omega}}}\left[  \left(\sum_{i=1}^{m_{\Omega}}\epsilon_i\Delta\phi_k(x_i)\right)^2 \right]\right)^{\frac{1}{2}}.
\end{align*}
where the first inequality is obtained due to the boundedness of the output weights in $2$ norm. Noticing that $\mathbb{E}[\epsilon_i\epsilon_j]=\delta_{ij}$, we further have
\begin{align*}
\mathbb{E}_{\epsilon_1,\cdots,\epsilon_{m_{\Omega}}}\mathbb{E}_{x_1,\cdots,x_{m_{\Omega}}}\left[ \left(\sum_{i=1}^{m_{\Omega}}\epsilon_i\Delta\phi_k(x_i)\right)^2\right] =\mathbb{E}_{x_1,\cdots,x_{m_{\Omega}}}\left[ \mathbb{E}_{\epsilon_1,\cdots,\epsilon_{m_{\Omega}}}\left(\sum_{i=1}^{m_{\Omega}}\epsilon_i\Delta\phi_k(x_i)\right)^2 \right] = \mathbb{E}_{x_1,\cdots,x_{m_{\Omega}}}\left[ m_{\Omega} \Delta\phi^2_k(x_1)\right]\leq m_{\Omega}M^2,
\end{align*}
from which we obtain $R_{m_{\Omega}}\left(\Delta\mathcal{H}_m\right)\lesssim \sqrt{mM^2}$. Similarly, we have $R_{m_{\Omega}}\left(\lambda f\circ\mathcal{H}_m\right)\lesssim \sqrt{m}$.
Therefore, we obtain the conclusion noticing that $M=O(m^\frac{3}{n+1+2\eta})$.
\end{proof}

Combining all of previous results, we obtain the following error estimate.
\begin{thm}
With the previous assumptions, we have the following local error estimate
\[ \| u^*- u^*_{md} \|^2_{2}  \leq C  m^{-\frac{2\eta}{n+2\eta +1}} +C_mm_p^{-\frac{1}{2}} \]
for some constant $C$ and $C_m=O(m^{\frac{3}{2}+\frac{3}{n+1+2\eta}})$, where $m$ is the number of neurons in the hidden layer, and the error estimate is understood in mean.
\end{thm}
\begin{remark}
So to guarantee the convergence, at least we need the number of sampling points growing far more fast than the number of neurons in the hidden layer. However, the rate of $C_m$ depending on $m$ can be further decreased by using different tools. Here we just present the worst case.
\end{remark}

\section{Improvements using the multilevel method}
In this section, in order to accelerate the calculation as well as improving the solution accuracy, adapting the multilevel method proposed in \cite{ying2024an}, we further propose a projection method by considering a minimization problem to fulfill the purpose of increasing the robustness of the algorithm.

For the discretized characteristic matrix $\Phi=[\Phi_{ij}]_{m_{\Omega}\times m}$, where $m_{\Omega}$ denotes the number of sampling points and $m$ is the number of neurons in the hidden layer (i.e., the number of basis functions), to carry out the initialization process followed by the optimization algorithm, relative large values of both $m$ and $m_{\Omega}$ would be needed for finding more solutions with better accuracy. However, increasing these two numbers too much would significantly increase the memory consumptions and the computational costs. Meanwhile, for the given threshold $R_{\epsilon}$ to determine accepting one approximation as a solution or not, on a coarse neural network, it is highly possible that an approximation is denied due to its low accuracy. So in order to partially overcome these difficulties, we firstly propose a projection method and then combine with the multilevel method, which turns out the new method can sharply reduce the computational costs while preserving high solution accuracy, even being able to find more solutions as shown in the numerical tests.

To be specific, given relatively small numbers of both $m^1$ and $m^1_{\Omega}$, by the initialization algorithm and the Newton method proposed in the previous sections, assume we have obtained one solution $\zeta^1$. Now with relatively large numbers of $m^2$ and $m^2_{\Omega}$, given a $\zeta^1$ and to map the solution from a coarse network to a fine network,  we firstly consider the following linear least square problem
\begin{equation}
\label{min-c-f}
\min_{\zeta^2\in\mathbb{R}^{m_2}}\frac{1}{N}\| H^1\zeta^1 - H^2\zeta^2 \|_{d}^2,
\end{equation}
where $H^1=\left(\phi^1_1, \phi^1_2, \cdots, \phi^1_{m_1}\right)$ and $H^2=\left(\phi^2_1, \phi^2_2, \cdots, \phi^2_{m_2}\right)$, the row vectors consisting of the randomized basis functions for the neural networks with the width $m_1$ and $m_2$, respectively,  $N$ is number of (uniformly distributed) sampling points to do the fitting, and the corresponding norm is defined to be
$\|u\|^2_{d} = \sum_{k=1}^N u^2(x_k)$ for $\forall  u\in C^{\infty}(\Omega)$. 
For simplicity, we can take $N=m^2_{\Omega}$ and use the sampling points in the fine network to calibrate the closeness of the solutions in different levels of networks. Like before, for this linear square problem, the solution can be easily obtained by solving the corresponding normal equation. Likewise, given a $\zeta^2$, we could also define another similar linear least square problem as follows
\begin{equation}
\label{min-f-c}
\min_{\zeta^1\in\mathbb{R}^{m_1}}\frac{1}{m^1_{\Omega}}\| H^1\zeta^1 - H^2\zeta^2 \|_{d}^2,
\end{equation}
which maps the solution from a fine network to a coarse network, as we will use later.
 Inspired by the Gauss-Newton method, suppose the solution of the least square problem \eqref{min-c-f} is denoted by $\zeta^2_{(0)}$, then as a good initial guess, we have

\[f(x, u(x)) \approx f(x, H^2\zeta^2_{(0)}(x)) + f'(x,H^2\zeta^2_{(0)}(x))(u(x) - H^2\zeta^2_{(0)}(x)),\]
from which we obtain
$u(x) \approx H^2\zeta^2_{(0)}(x) + H^2\zeta (x)$, where $\zeta$ satisfies
\begin{equation}\nonumber
\left\{\begin{array}{cl}
  \left[-(\Delta H^2) \zeta  - \lambda f'(x_j,H^2\zeta^2_{(0)}(x_j))H^2\zeta \right](x_j)=\left[\Delta H^2 \zeta^2_{(0)}+ \lambda f(x,H^2\zeta^2_{(0)})\right](x_j),~~&~~\{x_j\}_{j=1}^{m^2_{\Omega}}\in\Omega,\\
  H^2\zeta (x_j) = 0,~~&~\{x_j\}_{k=1}^{m^2_{\partial\Omega}}\in\partial\Omega,
\end{array}\right.
\end{equation}
where the equality is understood also in the least square sense. If this works, then we can repeat several times to further improve the solution's accuracy. Inspired by the multigrid method, we define a new method, called hierarchical randomized neural networks, to efficiently and accurately find multiple solutions of the considered semilinear problem. We describe the method in Algorithm 3 for the two level version, which is very like the full multigrid method. 

In the algorithm, for the initial solution, we can just take the trivial solution if it is. Otherwise, with an arbitrary initial guess, we need solve a nonlinear least square problem once using a randomized neural network. This process can be repeated to compute multiple solutions on a very fine network. And also due to its ability of improving accuracy,  we can set a relatively large threshold $R_{\epsilon}$ on the first course network and use a relatively small one on the fine networks to identify more potential solutions. 

	\begin{center}
		\begin{minipage}{0.98\linewidth}
			\hrule
			\vspace{0.5em}
			\noindent\textbf{Algorithm 3: Hierarchical randomized neural networks (Two levels) }
			\vspace{0.5em}
			
			\noindent\textbf{Input:}
			Two randomized neural networks with $m^1$ and $m^2$ hidden layers ($m^1\leq m^2$), $m^1_{\Omega}$ and $m^2_{\Omega}$ sampling points ($m^1_{\Omega}\leq m^2_{\Omega}$), respectively, and a given constant $R_{\epsilon}$. And the set of finding solutions is denoted by $U_s$. Suppose an initial solution has been identified, denoted by $\zeta_0$, and $U_s=\{\zeta_0\}$.
			
			\begin{enumerate}
			        \item For currently new solution or just added one in $U_s$, perform the initialization in Algorithm 1 and then the optimization in Algorithm 2 on the coarse network with the width $m^1$ using $m^1_{\Omega}$ sampling points. Here assume we find some approximations, the set of which is denoted by $\mathcal{A}$.
				\item In order to enhance the solution's accuracy, for each potential solution $\zeta^1\in \mathcal{A}$,
				\begin{enumerate}
				\item we do the prolongation operator by solving the linear least square problem \eqref{min-c-f} and obtain an improved approximation $\zeta^2$ by the Gauss-Newton method on the fine network with the width $m^2$ using $m^2_{\Omega}$ sampling points.
				\item Accept the improved approximation $\zeta^2$ as the solution if the norm of residues is less than or equal to $R_{\epsilon}$, and set $U_s=U_s\cup\{\zeta^2\}$; Otherwise deny the approximation as one solution.
				\end{enumerate}
				\item For each new solution $\zeta^2$, we do the restriction operator by solving the linear least square problem \eqref{min-f-c} on the coarse network with the width $m^1$ using $m^1_{\Omega}$ sampling points. And then back to step 1 to repeat the process to find more solutions .
			\end{enumerate}
			
			\noindent\textbf{Output:}
			The multiple solutions in $U_s$ with satisfied accuracy.
			\vspace{0.5em}
			\hrule
		\end{minipage}
	\end{center}

\section{Numerical experiments}

In this section, we carry out the numerical experiments using the proposed method. Here we use $u(x) = u(x_1), u(x)=u(x_1,x_2)$ and $u(x)=u(x_1,x_2,x_3)$ when $ x\in\mathbb{R}, x\in\mathbb{R}^2$ and $x\in\mathbb{R}^3$, respectively. All
the experiments are implemented by Python on a desktop PC with an Intel Xeon E5-2666 V3 CPU (2.90 GHz) and 32 GB of RAM. The shallow neural networks use $\sigma = tanh(x+0.5)-tanh(x-0.5)$ as the activation function. And the weights $w_1$ and bias $b_1$ of the neural networks are initialized by the scrambled Sobol sequences. That is, $
\mathbf{w}_1=s(2\mathbf{Q}-\mathbf{1})$,
$\mathbf{b}_1=-(\mathbf{w}_1\odot\boldsymbol{\Xi})\mathbf{1},
$
where $\mathbf{Q}$ and $\boldsymbol{\Xi}$ denote the two coordinate blocks of a scrambled Sobol point set, $\odot$ denotes the element-wise product, and $s$ is the weight scaling parameter. In the coarse level, we set $s=4$, and in the refinement levels, we set $s=8$ to improve the expressivity of the neural networks. Unless otherwise specified, the structures of the shallow neural networks as well as the sampling points in different levels follow the numbers given in  Table~\ref{tab:network_config}.

\begin{table}[!htb]
	\caption{Numbers of random features and collocation points.}
	\label{tab:network_config}
	\centering
	\normalsize
	\setlength{\tabcolsep}{9pt}
	\renewcommand{\arraystretch}{1.2}
	\begin{tabular}{lrrrrrr}
		\toprule
		\rowcolor{white}
		& \multicolumn{2}{c}{1D}
		& \multicolumn{2}{c}{2D}
		& \multicolumn{2}{c}{3D} \\
		\cmidrule(lr){2-3}
		\cmidrule(lr){4-5}
		\cmidrule(lr){6-7}
		Level
		& $M$ & $M_\Omega$
		& $M$ & $M_\Omega$
		& $M$ & $M_\Omega$ \\
		\midrule
		Coarse  & 120 & 400
		& 600 & 2000
		& 720 & 2400 \\
		Level 1 & 240 & 1200
		& 1000 & 4000
		& 1300 & 4800 \\
		Level 2 & 400 & 2400
		& 1600 & 8000
		& 2200 & 8000 \\
		Level 3 & 640 & 4000
		& 3500 & 12000
		& 3500 & 12000 \\
		\bottomrule
	\end{tabular}
\end{table}

In the process of initializations, the parameters $\epsilon_1$ is set as $10^{-11}$ at the coarse level and $5\times10^{-16}$ at the refinement levels. In all levels, we employ $K=6$ and $\alpha_{init}=0.5, 1, 2, 4, 8$. The deflation parameters are set as $\theta =2$ and $\sigma = 1$. In GDGN, we set $R_{\epsilon}=0.1$ and employ a line search with an initial step size of \(\alpha=1\), and the step size is halved at most 16 times until a reduction in the residual is achieved. Meanwhile, in order to test the accuracy of the proposed method, we use the following relative $L^2$ errors for problems with exact solutions available and RMS errors for others,

\begin{equation}
	\mathrm{Err}_{L^2}
	=
	\frac{
		\sqrt{\sum_{i=1}^{m}
			|\hat{u}(\boldsymbol{x}_i)-u(\boldsymbol{x}_i)|^2}
	}{
		\sqrt{\sum_{i=1}^{m}
			|u(\boldsymbol{x}_i)|^2}
	},\qquad
	\mathrm{Err}_{\mathrm{RMS}}
	=
	\left(
	\frac{1}{m}
	\sum_{i=1}^{m}
	|\mathcal{F}[\hat{u}](\boldsymbol{x}_i)|^2
	\right)^{1/2},
\end{equation}
where $\hat{u}$ and $u$ represent the numerical solutions and the exact solutions, respectively, 
$\mathcal{F}[\hat{u}]$ denotes the residual of the PDEs, and $m$ represents the total number of test points chosen uniformly in the domain. For the 1D, 2D and 3D problems, we set $m=2^{12}$, $2^{13}$, and $2^{14}$,  respectively. 

\begin{example}[Bratu problem]
	We first consider the following Bratu problem, 
	\begin{equation}\label{Example1}
		\begin{cases}
			-u''=\lambda e^{u}, & x\in \Omega,\\
			u=0, & x \in \partial\Omega,
		\end{cases}
	\end{equation}
	where $\lambda=0.5$ and $\Omega=(-1, 1)$. This problem has two exact solutions given by $
		u=
		\log\left(
		z^2\operatorname{sech}^2
		\left(\frac{zx}{2}\right)
		\right)$, 
	where $z$ takes either of the two positive roots of
	$z=\cosh\left(\frac{z}{2}\right)$.

	\begin{table}[!htb]
		\centering
		\caption{Numerical errors of the multilevel solutions for Example 1.}
		\label{tab:ex1_multilevel}
		\small
		\renewcommand{\arraystretch}{1.1}
		\setlength{\tabcolsep}{5pt}
		\begin{tabular}{lccccc}
			\toprule
		        \rowcolor{white}
			\diagbox[width=7em,height=2.2em]{Solution}{Level}
			& Coarse & Level 1 & Level 2 & Level 3
			& $Err_{L^2}$ \\
			\midrule
			Solution 1 & 3.718e-05 & 1.276e-13 & 9.410e-14
			& 6.097e-14 & 1.570e-15 \\
			Solution 2 & 1.545e-02 & 4.579e-13 & 2.643e-13
			& 2.093e-13 & 4.060e-14 \\
			\bottomrule
		\end{tabular}
	\end{table}

	Table~\ref{tab:ex1_multilevel} reports the training errors on different levels. It's clearly observed that the multilevel refinement significantly reduces the PDE residuals, thus yielding the final relative \(L^2\) errors of \(1.57\times10^{-15}\) and \(4.06\times10^{-14}\), respectively. These results demonstrate the effectiveness of the proposed method in identifying multiple solutions and achieving high numerical accuracy. For completeness, the numerical solutions and the exact solutions as well as the residuals are presented in Fig~\ref{Fig.1}.

	\begin{figure}[!htb]
		\centering
		\includegraphics[width=0.7\linewidth]{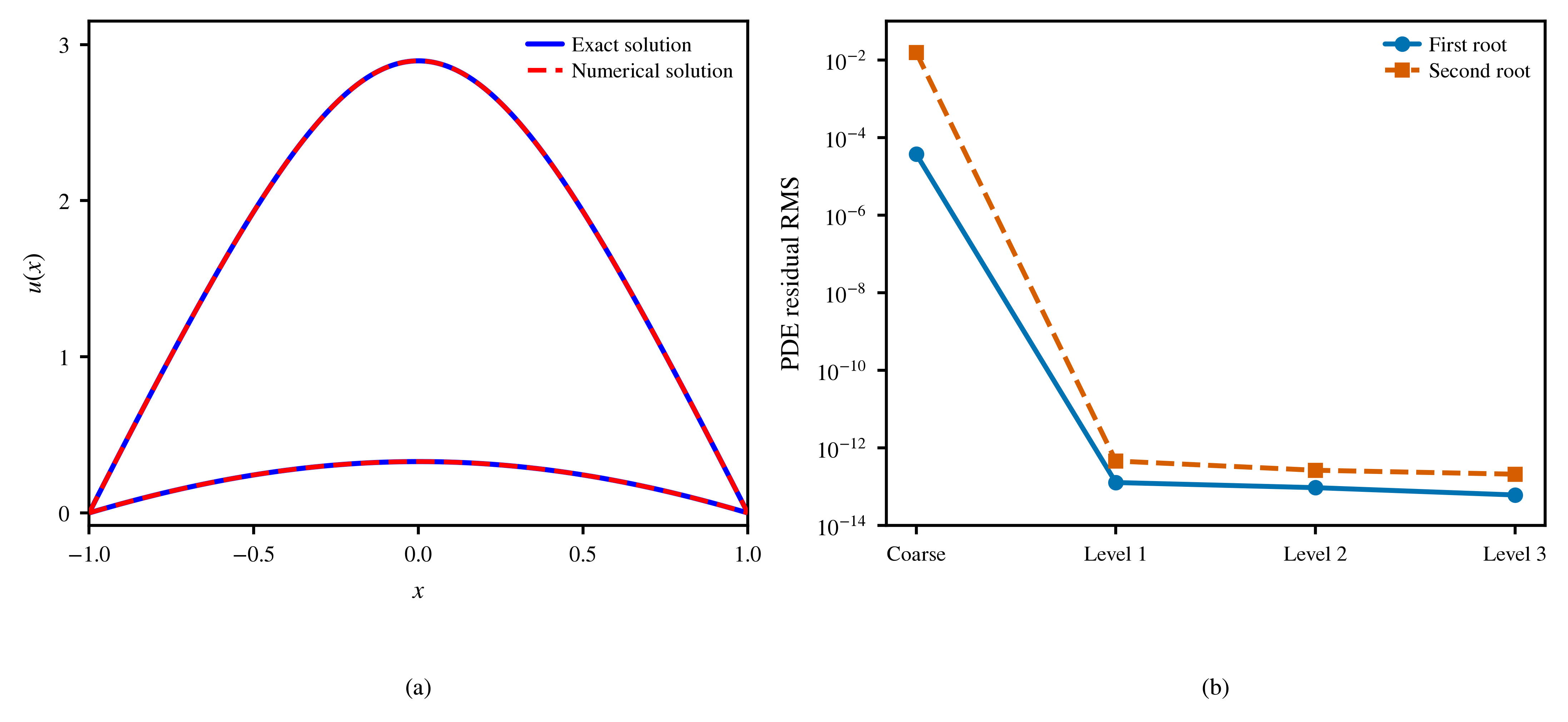}
		\caption{Numerical results for Example 1.
			(a) Comparison of numerical and exact solutions;
			(b) PDE residuals in different levels.}
		\label{Fig.1}
	\end{figure}

\end{example}

\begin{example} Next we consider the following boundary value problem, 
	\begin{equation}\label{Example2}
		\begin{cases}
			(u'')^2-(1+e^x)u''+e^x=0,
			& x\in \Omega,\\
			u=0, & x \in \partial \Omega.
		\end{cases}
	\end{equation}
	where $\Omega = (0, 1)$, and exact solutions are 
			$u_1=\frac{x^2-x}{2}$ and 
			$u_2=e^x-1+(1-e)x$.

	\begin{table}[!htb]
		\centering
		\caption{Numerical errors of the multilevel solutions for Example 2.}
		\label{tab:example2}
		\small
		\renewcommand{\arraystretch}{1.1}
		\setlength{\tabcolsep}{5pt}
		\begin{tabular}{lccccc}
			\toprule
			\rowcolor{white}
			\diagbox[width=7em,height=2.2em]{Solution}{Level}
			& Coarse & Level 1 & Level 2 & Level 3
			& \shortstack{$Err_{L^2}$} \\
			\midrule
			Solution 1  & 1.615e-03 & 2.924e-13 & 2.275e-13
			& 4.666e-13 & 2.210e-15 \\
			Solution 2  & 4.429e-06 & 2.174e-12 & 1.768e-12
			& 3.304e-12 & 2.060e-15 \\
			\bottomrule
		\end{tabular}
	\end{table}
	Table~\ref{tab:example2} reports the training PDE residuals at different refinement levels for Example 2. We can see that the proposed method successfully identifies both solutions, with relative \(L^2\) errors on the order of \(10^{-15}\). Although the residuals increase slightly at Level 3, they remain at the order of \(10^{-12}\)–\(10^{-13}\). This may be caused by roundoff errors, as the residuals have already reached a very low level.
	
\end{example}

\begin{example}
	\label{Example3}
	Here we consider the following nonlinear boundary value problem,
	\begin{equation}
		\begin{cases}
			(u'')^2+u''-2=0, & x\in \Omega,\\
			u=0, & x \in \partial\Omega,
		\end{cases}
	\end{equation}
	where $\Omega=(0,1)$. This nonlinear equation can be factorized as $(u''-1)(u''+2)=0$, which has two distinct
	classical solutions under the homogeneous
	Dirichlet boundary conditions, where
			$u_1=\frac{x^2-x}{2}$ and $u_2=x(1-x)$.
	This example also provides a simple test of the proposed method's ability to identify distinct solution branches.
	
	From Table~\ref{tab:example3}, we can see the training PDE residuals at different refinement levels for Example 3. In the coarse level, the residuals of two solutions are $2.329\times 10^{-9}$ and $2.734\times 10^{-3}$. After the first refinement, both residuals decrease to $10^{-12}$, indicating a substantial improvement. Like the previous example, the residuals increase slightly at Level 3, which may be caused by the roundoff errors, as the solutions have already achieved very high accuracy. Moreover, the relative $L^2$ errors reach $10^{-16}$, approaching the machine precision. All of these results show the accuracy of the proposed multilevel method.

	\begin{table}[!htb]
		\centering
		\caption{Numerical errors of the multilevel solutions for Example 3.}
		\label{tab:example3}
		\normalsize
		\renewcommand{\arraystretch}{1.15}
		\setlength{\tabcolsep}{6pt}
		\begin{tabular}{lccccc}
			\toprule
			\rowcolor{white}
			\diagbox[width=7em,height=2.5em]{Solution}{Level}
			& Coarse & Level 1 & Level 2 & Level 3
			& \shortstack{$Err_{L^2}$} \\
			\midrule
			Solution 1 & 2.329e-09 & 3.830e-12 & 1.441e-12
			& 2.011e-12 & 4.240e-16 \\
			Solution 2 & 2.734e-03 & 7.660e-12 & 2.882e-12
			& 4.022e-12 & 4.100e-16 \\
			\bottomrule
		\end{tabular}
	\end{table}
	
\end{example}

\begin{example}
	\label{Example4}
	Given a non-negative parameter $p$, We consider the following nonlinear boundary value problem 
	\begin{equation}
		\begin{cases}
			u''=u^2(u^2-p), & x\in(0,1),\\
			u'(0)=0,\quad u(1)=0,
		\end{cases}
	\end{equation}	
	Here in order to enforce the mixed boundary conditions, we employ
	the modified random features defined by
	$\phi(x)=(1-x^2)\sigma(w_1x^2+b_1)$. Meanwhile, unlike the previous examples, we set $R_{\epsilon}=10$, the coarse-level feature scale $s = 20$, and the refinement feature scales to $(32,36,40)$. 
	
	For this test case, as pointed out in \cite{huang2022hompinns}, for any given parameter $p$, there exist multiple solutions, and the number of solutions increases as $p$ goes large. So here we consider $p=3,10,18$, respectively, to see the abilities of our method identifying the different solutions. As shown in Fig~\ref{fig:ex4_solutions},
	the proposed method successfully identifies all nontrivial solutions for each $p$, demonstrating its effectiveness in discovering multiple solutions.
	\begin{figure}[!htb]
		\centering
		\includegraphics[width=\linewidth]{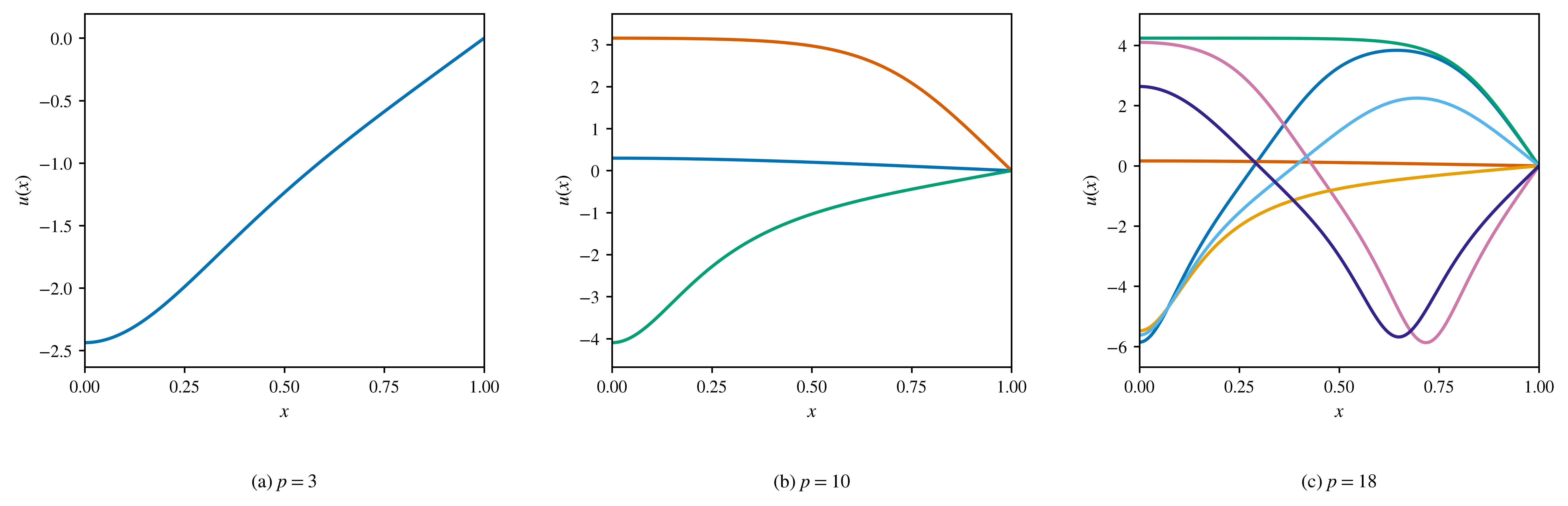}
		\caption{Nontrivial numerical solutions of Example 4 for different $p$.}
		\label{fig:ex4_solutions}
	\end{figure}
	
	Table~\ref{tab:example4} presents the numerical results for $p=10$. The proposed method successfully identifies all three nontrivial solutions. Multilevel refinement significantly reduces the PDE residuals, with two solutions reaching the order of $10^{-10}$--$10^{-12}$, while the third achieves $10^{-4}$.

	\begin{table}[!htb]
		\centering
		\caption{Numerical errors of the nontrivial multilevel solutions for Example 4 with $p=10$.}
		\label{tab:example4}
		\normalsize
		\renewcommand{\arraystretch}{1.15}
		\setlength{\tabcolsep}{6pt}
		\begin{tabular}{lccccc}
			\toprule
			\rowcolor{white}
			\diagbox[width=6.5em,height=2.5em]{Solution}{Level}
			& Coarse & Level 1 & Level 2 & Level 3
			& \shortstack{$Err_{RMS}$} \\
			\midrule
			Solution 1 & 1.702e-03 & 9.768e-11 & 3.945e-10
			& 1.980e-10 & 2.031e-10 \\
			Solution 2 & 4.030e-02 & 4.059e-12 & 1.883e-11
			& 9.506e-12 & 9.737e-12 \\
			Solution 3 & 6.708e-02 & 9.573e-04 & 3.323e-04
			& 2.999e-04 & 3.014e-04 \\
			\bottomrule
		\end{tabular}
	\end{table}
	
\end{example}

\begin{example}
	\label{Example5}
	In this case, we consider the following 2D nonlinear elliptic problem,
	\begin{equation}
		\begin{cases}
			-\Delta u=u^3, & x\in\Omega,\\
			u=0, & x\in\partial\Omega,
		\end{cases}
	\end{equation}
	where $\Omega=(0,1)^2$. This problem has a trivial solution $u\equiv 0$ and infinitely many non-trivial solutions. Meanwhile, the set of the exact solutions is invariant under the dihedral group $D_4$, so that the rotation of the solution by integer multiple of $90$ degree is still a solution. 
	
	In this case, we set coarse-level threshold $R_{\epsilon}=5$ to make more candidates for refinement. As shown in Table~\ref{tab:example5}, although some coarse solutions have relatively large residuals, multilevel refinement reduces them by several orders of magnitude.	
	\begin{table}[!htb]
		\centering
		\caption{Numerical errors of the nontrivial multilevel solutions for Example 5.}
		\label{tab:example5}
		\normalsize
		\renewcommand{\arraystretch}{1.15}
		\setlength{\tabcolsep}{6pt}
		\begin{tabular}{lccccc}
			\toprule
			\rowcolor{white}
			\diagbox[width=6.6em,height=2.5em]{Solution}{Level}
			& Coarse & Level 1 & Level 2 & Level 3
			& \shortstack{$Err_{RMS}$} \\
			\midrule
			Solution 1  & 1.149e+00 & 3.831e-07 & 9.237e-09
			& 2.127e-10 & 2.632e-10 \\
			Solution 2  & 4.675e+00 & 3.831e-07 & 9.237e-09
			& 2.127e-10 & 2.632e-10 \\
			Solution 3  & 4.457e+00 & 3.291e-03 & 2.313e-04
			& 1.038e-06 & 1.283e-06 \\
			Solution 4  & 4.467e+00 & 3.291e-03 & 2.313e-04
			& 1.038e-06 & 1.283e-06 \\
			Solution 5  & 4.565e+00 & 4.054e-03 & 1.843e-04
			& 1.094e-06 & 1.264e-06 \\
			Solution 6  & 4.575e+00 & 4.054e-03 & 1.843e-04
			& 1.094e-06 & 1.264e-06 \\
			Solution 7  & 4.869e+00 & 5.839e-02 & 5.252e-03
			& 4.932e-05 & 5.758e-05 \\
			Solution 8  & 3.003e+00 & 5.839e-02 & 5.252e-03
			& 4.932e-05 & 5.758e-05 \\
			Solution 9  & 4.189e+00 & 1.020e-01 & 1.048e-02
			& 1.362e-04 & 1.566e-04 \\
			Solution 10 & 4.399e+00 & 1.174e-01 & 1.163e-02
			& 1.434e-04 & 1.579e-04 \\
			Solution 11 & 4.144e+00 & 1.174e-01 & 1.163e-02
			& 1.434e-04 & 1.579e-04 \\
			Solution 12 & 4.486e+00 & 1.020e-01 & 1.048e-02
			& 1.362e-04 & 1.566e-04 \\
			\bottomrule
		\end{tabular}
	\end{table}
	Meanwhile, to further show the effectiveness of the multilevel method, we compare it with direct nonlinear optimization at Level 2 (Here the direct nonlinear optimization at Level 3 would be computational costs and also may exceed the memory limit). So in order to make a fair comparison, all other procedures remain unchanged, including initialization and solution search algorithms.
	
	Table~\ref{tab:ex7_multilevel_comparison} presents the results, which demonstrates that the proposed multilevel method not only increases the solution accuracy significantly using less than one tenth computational cost,  but also is able to find more distinct solutions. This is particular important when we search for as many solution as possible. 	
	\begin{table}[!htb]
		\centering
		\caption{Comparison results for Example 5.}
		\renewcommand{\arraystretch}{1.15}
		\label{tab:ex7_multilevel_comparison}
		\setlength{\tabcolsep}{6pt}
		\begin{tabular}{lcccc}
			\toprule
			\rowcolor{white}
			Method & Solutions & Min. $\mathrm{Err}_{\mathrm{RMS}}$
			& Max. $\mathrm{Err}_{\mathrm{RMS}}$ & Time (s) \\
			\midrule
			Direct Level-2
			& 8 & $2.803\times10^{-4}$ & $6.141\times10^{-3}$ & 94121 \\
			Multilevel
			& 12 & $2.724\times10^{-10}$ & $1.659\times10^{-4}$ & 7093 \\
			\bottomrule
		\end{tabular}
	\end{table}
	Lastly, we also show the structures of the identified nontrivial solutions in Fig~\ref{fig:ex5_solutions}, from which we can see that all symmetric solutions can be identified. That is, as long as $u(x_1,x_2)$ is a solution, so is the $D_4u(x_1,x_2)$ for the dihedral group $D_4$. This is another evaluable point demonstrating the effectiveness of our proposed method.

	\begin{figure}[!htb]
		\centering
		\includegraphics[width=0.8\linewidth]{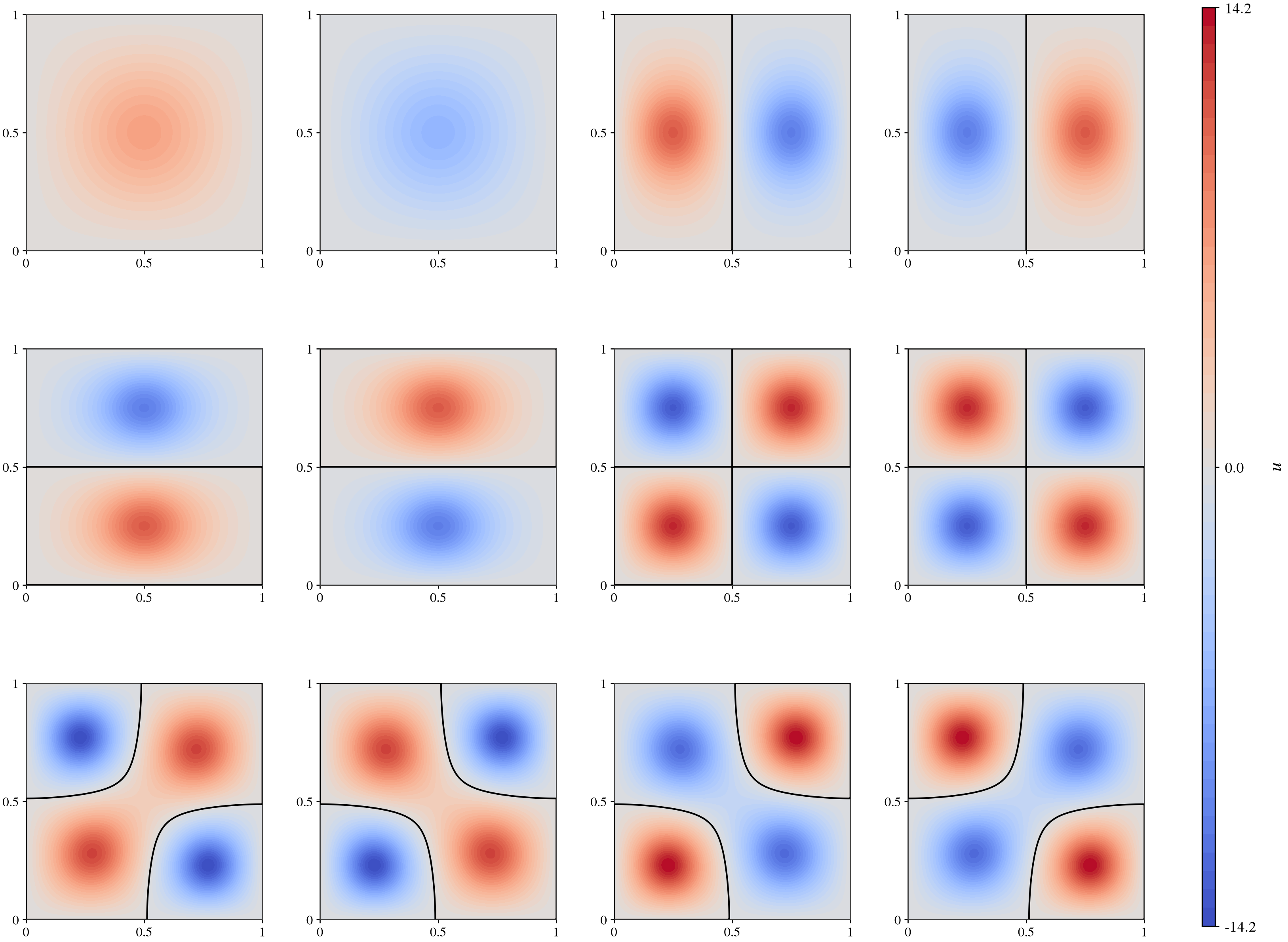}
		\caption{Nontrivial numerical solutions of Example 5.}
		\label{fig:ex5_solutions}
	\end{figure}
	
\end{example}

\begin{example}
	\label{Example7}
	We consider the following two-dimensional
	Lotka-Volterra competition system,
	\begin{equation}
		\begin{cases}
			-\Delta u=u(30-u-0.2v), & x\in\Omega,\\
			-\Delta v=v(30-0.2u-v), & x\in\Omega,\\
			u=v=0, & x\in\partial\Omega,
		\end{cases}
	\end{equation}
	where $\Omega=(0,1)^2$.
	This problem has four nonnegative solutions, including a trivial solution, two semitrivial solutions, and one positive coexistence solution. Each component uses the same number of random features as in the previous 2D examples, and the two equations are solved simultaneously.
	
	As shown in Table~\ref{tab:example6}, the proposed method again successfully identifies all four nonnegative solutions. The multilevel refinement again reduces the PDE residuals to the order of $10^{-8}$, with $Err_{RMS}$ of order $10^{-7}$. All of these results show the effectiveness of our proposed method again for the system of nonlinear problems.

	\begin{table}[!htb]
		\centering
		\caption{Numerical errors of the nontrivial multilevel solutions for Example 6.}
		\label{tab:example6}
		\normalsize
		\renewcommand{\arraystretch}{1.15}
		\setlength{\tabcolsep}{6pt}
		\begin{tabular}{lccccc}
			\toprule
			\rowcolor{white}
			\diagbox[width=6.5em,height=2.5em]{Solution}{Level}
			& Coarse & Level 1 & Level 2 & Level 3
			& \shortstack{$Err_{RMS}$} \\
			\midrule
			Solution 1 & 1.505e-03 & 9.159e-07 & 5.014e-07
			& 7.607e-08 & 1.188e-07 \\
			Solution 2 & 2.817e-04 & 9.159e-07 & 5.014e-07
			& 7.607e-08 & 1.188e-07 \\
			Solution 3 & 3.017e-03 & 1.079e-06 & 5.910e-07
			& 8.965e-08 & 1.400e-07 \\
			\bottomrule
		\end{tabular}
	\end{table}
	
\end{example}

\begin{example}
	\label{Example6}
	Here we consider the following two-dimensional stationary
	Allen--Cahn equation,
	\begin{equation}
		\begin{cases}
			-\varepsilon \Delta u + u^3-u=0,
			& (x,y)\in (0, 1)^2,\\
			u(x,0)=u(x,1)=-4x(1-x), & x\in[0,1]\\
			u(0,y)=u(1,y)=4y(1-y), & y\in [0,1]
		\end{cases}
	\end{equation}
	where $\varepsilon=0.01$. Different from the previous examples, this problem is subject to nonhomogeneous Dirichlet boundary conditions. So in this case, the soft constraint is used to handle the nonhomogeneous boundary conditions.
	
	From Fig~\ref{fig:ex7_solutions}, we can see that
	our method successfully identifies three distinct solutions. Table~\ref{tab:example7} presents the PDE residuals at different refinement levels. The results shows that without requiring any prior knowledge of the solution structures, our method once again successfully discovers multiple solutions and achieves high accuracy. The validation RMS errors of the three solutions are $3.892\times10^{-8}$, $4.090\times10^{-8}$,
	and $3.981\times10^{-8}$, respectively. These results further demonstrate the effectiveness and accuracy of the proposed method in solving nonlinear PDEs with multiple solutions.

	\begin{table}[!htb]
		\centering
		\caption{Numerical errors of the multilevel solutions
			for Example 7.}
		\label{tab:example7}
		\normalsize
		\renewcommand{\arraystretch}{1.15}
		\setlength{\tabcolsep}{6pt}
		\begin{tabular}{lccccc}
			\toprule
			\rowcolor{white}
			\diagbox[width=6.5em,height=2.5em]{Solution}{Level}
			& Coarse & Level 1 & Level 2 & Level 3
			& $Err_{RMS}$ \\
			\midrule
			Solution 1
			& 7.938e-03 & 5.176e-07 & 1.166e-07
			& 2.535e-08 & 3.892e-08 \\
			Solution 2
			& 3.812e-02 & 3.044e-06 & 2.553e-07
			& 2.580e-08 & 4.090e-08 \\
			Solution 3
			& 8.285e-03 & 3.717e-06 & 2.085e-07
			& 2.600e-08 & 3.981e-08 \\
			\bottomrule
		\end{tabular}
	\end{table}

	\begin{figure}[!htb]
		\centering
		\includegraphics[width=0.7\linewidth]{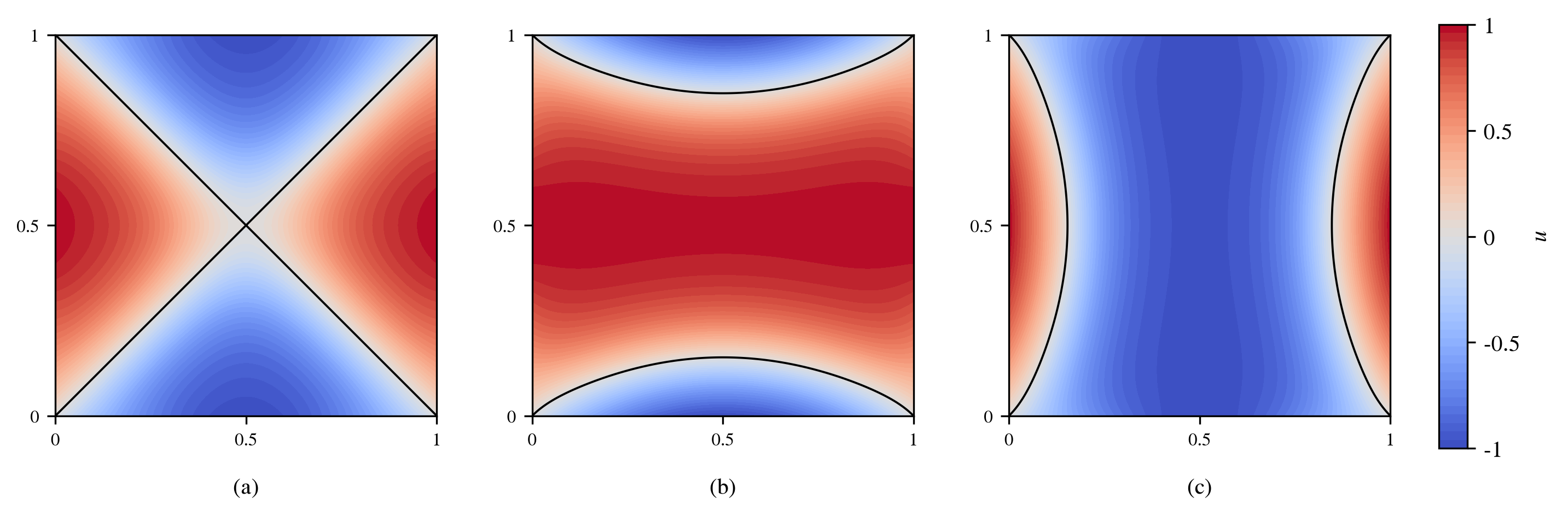}
		\caption{Results of Example 7. Here (a) shows the unstable solution, while (b) and (c) show the positive and negative stable solutions, respectively.}
		\label{fig:ex7_solutions}
	\end{figure}

\end{example}

\begin{example}
	\label{Example8}
	We lastly consider the following three-dimensional Allen-Cahn equation
	\begin{equation}
		\begin{cases}
			-\Delta u=\lambda u-u^3, & x\in\Omega,\\
			u=0, & x\in\partial\Omega,
		\end{cases}
	\end{equation}
	where $\Omega=(0,\pi)^3$ and $\lambda=5$. For $\lambda=5$, this problem has three solutions: the trivial solution $u=0$ and a pair of symmetric positive and negative solutions. 
	As shown in Fig~\ref{fig:ex8_solutions}, the proposed method successfully identifies all nontrivial solutions. And Table~\ref{tab:example8} presents the training loss and validation $Err_{RMS}$. From the table, it's clear that the multilevel refinement reduces the PDE residuals to $3.033\times10^{-5}$, with $Err_{RMS}$ of $3.704\times10^{-5}$ for both nontrivial solutions. Since these two nontrivial solutions are $u$ and $-u$, the symmetry known from the equation, and our method seems exactly preserve this hidden symmetry, which explaining why except the coarse level, it seems we have exactly the same residuals for each refined levels and also exact the same RMS for the obtained different solutions lastly.

	\begin{figure}[!htb]
		\centering
		\includegraphics[width=0.65\linewidth]{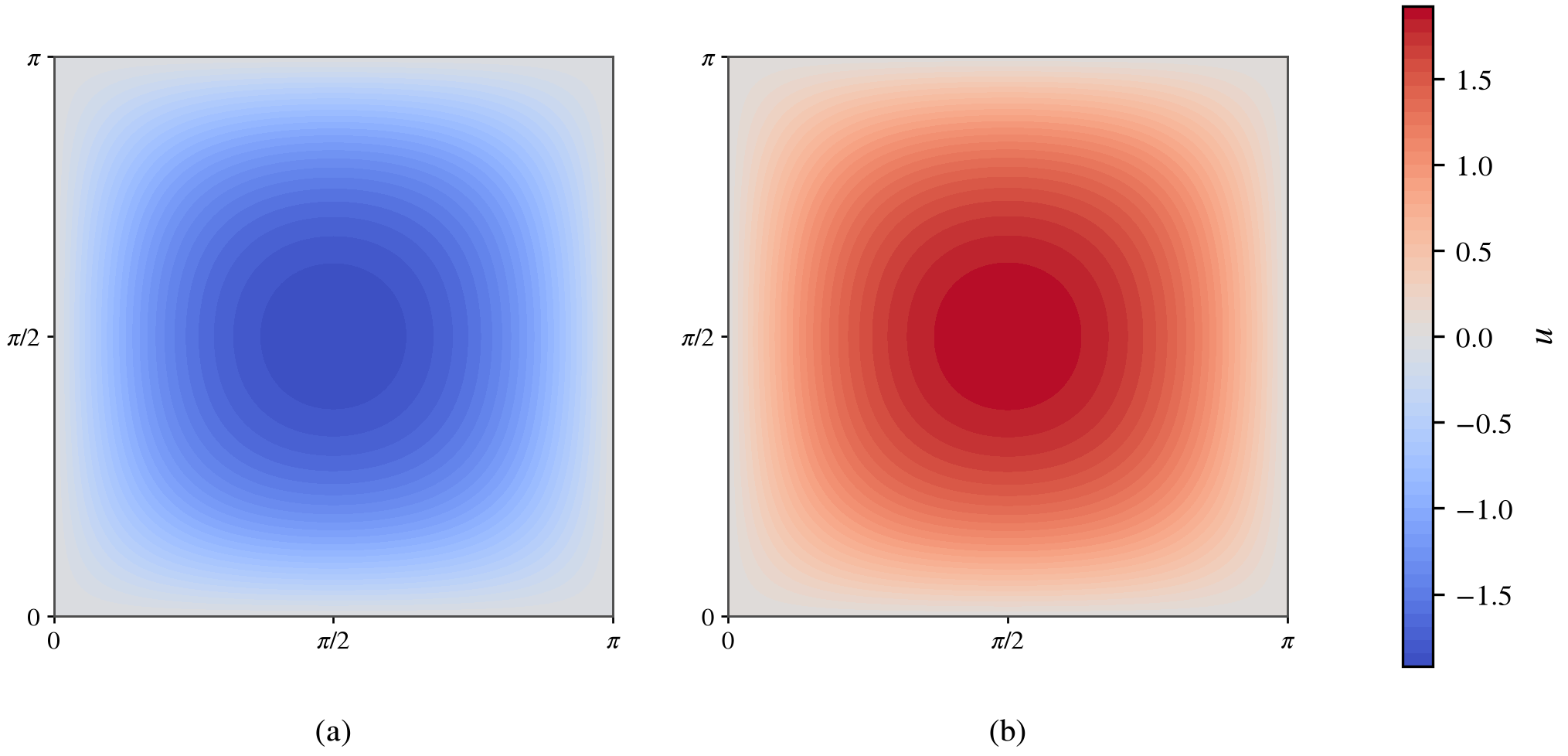}
		\caption{Nontrivial numerical solutions of Example 8 on the plane $z=\frac{\pi}{2}$. Here (a) and (b) show the negative and positive solutions, respectively.}
		\label{fig:ex8_solutions}
	\end{figure}

	\begin{table}[!htb]
		\centering
		\caption{Numerical errors of the nontrivial multilevel
			solutions for Example 8.}
		\label{tab:example8}
		\renewcommand{\arraystretch}{1.1}
		\setlength{\tabcolsep}{6pt}
		\begin{tabular}{lccccc}
			\toprule
			\rowcolor{white}
			\diagbox[width=7em,height=2.2em]{Solution}{Level}
			& Coarse & Level 1 & Level 2 & Level 3
			& $Err_{\mathrm{RMS}}$ \\
			\midrule
			Solution 1 & 2.088e-03 & 5.064e-04 & 1.107e-04
			& 3.033e-05 & 3.704e-05 \\
			Solution 2 & 6.069e-03 & 5.064e-04 & 1.107e-04
			& 3.033e-05 & 3.704e-05 \\
			\bottomrule
		\end{tabular}
	\end{table}
	
\end{example}

\section{Conclusion}
This paper was devoted to developing the method for efficiently and accurately finding multiple solutions of semilinear elliptic problems. The improvements were done mainly from three different aspects, initializations, optimization, and the multilevel technique for accelerations and enhancing solutions' accuracies. Local error estimate of a particular solution satisfying some conditions were also considered in terms of both numbers of the neurons in the hidden layer and the sampling points. Meanwhile, the numerical tests on several examples with or without exact solutions not only demonstrated the effectiveness of the proposed method (the solution accuracy reaching the machine accuracy for the best case, while the computational cost being reduced to less than $8\%$ for the tested case), but also preserved the inherent symmetries of multiple solutions, latter of which is significant in some applications. 

For the proposed method, high dimensional problems are not tested yet, which will be considered in the future. Meanwhile, how to fulfill the proposed algorithm in a cluster using some parallel techniques to identify multiple solutions of some large-scale engineering problems is also worth to consider. Lastly, as the first paper in a series of work, we will continue to work on semilinear problems with exponential nonlinearity or others, to see if further theoretical analysis on solutions can be used to alleviate the computation to find more solutions.
\section*{Acknowledgments}
The work was supported by the Natural Science Foundation of China (12671352), the Research Foundation and the Young Core Instructor Project of  Education Department of Hunan Province (25A0167, 25B0011, and HNJG-20230343), the Major Scientific and Technological Innovation Platform Project of Hunan Province (2024JC1003), and the Natural Science Foundation of Hunan Province (2025JJ50013). We are also grateful to the High Performance Computing Center of Central South University for partial support of this work.

\section*{Declaration of interest}
The authors declare that they have no known competing financial interests or personal relationships that could have appeared to influence the work.







\end{document}